\documentclass{article}
\usepackage{url} %
\usepackage{lineno}
\usepackage{amsmath,amsfonts,mathrsfs,amssymb,mathtools}
\usepackage{bbm}
\usepackage{tikz,pgfplots}
\usepackage{subfig}
\usepackage{graphbox}
\usepackage{float}
\usepackage{multirow}
\usepackage{booktabs}
\usepackage{amsthm}
\usepackage{algorithm}          %
\usepackage{float}
\newfloat{algorithm}{t}{lop}
\usepackage{geometry}
\usepackage{algpseudocode}      %
\newcommand{\bs}[1]{{\boldsymbol #1}}
\newcommand{\citeA}[1]{\cite{#1}}

\newcommand{\trans}{\mathsf{T}}
\newcommand{\transpose}{^{\trans}}

\newcommand{\gradient}{\nabla}

\newcommand{\divergence}{\nabla\cdot}

\newcommand{\norm}[1]{\left\lVert#1\right\rVert}

\newcommand{\boldvar}[1]{\ensuremath{\boldsymbol{#1}}}
\newcommand{\boldvec}[1]{\ensuremath{\mathbf{#1}}}
\newcommand{\mat}[1]{\boldvec{#1}}
\ifdefined\vec
  \renewcommand{\vec}[1]{\boldvec{#1}}
\else
  \newcommand{\vec}[1]{\boldvec{#1}}
\fi

\newcommand{\vel}{\boldvar{u}}
\newcommand{\press}{p}

\newcommand{\visc}{\mu}

\newcommand{\strainrate}{\dot{\varepsilon}}
\newcommand{\strainratetensor}{\dot{\boldvar{\varepsilon}}}

\newcommand{\viscstress}{\tau}
\newcommand{\viscstresstensor}{\boldvar{\tau}}

\newcommand{\invII}[1]{#1_\textsc{ii}} 
\newcommand{\strainrateinvII}{\invII{\strainrate}}
\newcommand{\viscstressinvII}{\invII{\viscstress}}

\newcommand{\rhsvel}{\boldvar{f}}

\newcommand{\viscref}{{\visc_{\mathrm{r}}}}
\newcommand{\viscmin}{\underline\visc}

\newcommand{\viscstressyield}{\viscstress_{\mathrm{y}}}

\definecolor{col-co}{RGB}{255,69,0} %
\definecolor{col-id}{RGB}{30,144,255} %

\begin{document}

\title{Advanced Newton methods for geodynamical models of Stokes flow with viscoplastic rheologies}

\author{%
Johann Rudi%
\thanks{Mathematics and Computer Science Division, Argonne National Laboratory, Lemont IL, USA (\texttt{jrudi@anl.gov})}%
\and
Yu-hsuan Shih%
\thanks{Courant Institute of Mathematical Sciences, New York University, New York NY, USA}%
\and
Georg Stadler%
\footnotemark[2]%
}
\maketitle

\begin{abstract}
\noindent
Strain localization and resulting plasticity and failure play an important role
in the evolution of the lithosphere.  These phenomena are commonly modeled by
Stokes flows with viscoplastic rheologies.  The nonlinearities of these
rheologies make the numerical solution of the resulting systems challenging,
and iterative methods often converge slowly or not at all.  Yet accurate
solutions are critical for representing the physics.  Moreover, for some
rheology laws, aspects of solvability are still unknown.  We study a basic but
representative viscoplastic rheology law.  The law involves a yield stress that
is independent of the dynamic pressure, referred to as von Mises yield
criterion.  Two commonly used variants, perfect/ideal and composite
viscoplasticity, are compared.  We derive both variants from energy
minimization principles, and we use this perspective to argue when solutions
are unique.  We propose a new stress--velocity Newton solution algorithm that
treats the stress as an independent variable during the Newton linearization
but requires solution only of Stokes systems that are of the usual
velocity--pressure form.
To study different solution algorithms, we implement 2D and 3D finite element
discretizations, and we generate Stokes problems with up to 7 orders of
magnitude viscosity contrasts, in which compression or tension results in
significant nonlinear localization effects.  Comparing the performance of the
proposed Newton method with the standard Newton method and the Picard
fixed-point method, we observe a significant reduction in the number of
iterations and improved stability with respect to problem nonlinearity, mesh
refinement, and the polynomial order of the discretization.
\end{abstract}

\section{Introduction}\label{sec:intro}
Numerical simulations are important for understanding models of lithosphere
dynamics over long time scales.
These models are often non-Newtonian and include nonlinear
mechanisms that can initiate %
norm faults \cite{LavierBuckPoliakov00, SharplesMoresiJadamecEtAl15},
strike-slip faults \cite{SobolevPetruninGarfunkelEtAl05}, and shear zones for
subduction \cite{NaliboffBillenTaras13, MaoGurnisMay17, Tackley00part1} through strain
localization. %
This is commonly achieved %
by incorporating
a form of frictional plasticity \cite{Ranalli95, Fossen16, Byerlee68}.
The transition from viscous to plastic
behavior is controlled by stress-limiting techniques using a yield stress,
whose value commonly depends on temperature and (lithostatic) pressure.
The resulting nonlinear models take the form of incompressible Stokes equations
and require numerical approximations of their solutions
that have to be accurate and reliable.
Obtaining accurate numerical %
solutions of models with plasticity is computationally challenging, however,
because of the
nonlinearities in the resulting systems
\cite{SpiegelmanMayWilson16, DuretzSoucheBorstEtAl18, PopovSobolev08, Kaus10}.
For instance, it can be difficult to
obtain converged solutions where the (nonlinear) residual of the system has
been reduced sufficiently.
Moreover, these solutions strongly depend on
the discretization, such as the refinement of the computational mesh.

In this paper, we study the solvability of viscoplastic Stokes flow in theory
and with numerical experiments. Our main contribution is the development of
nonlinear iterative methods with fast convergence that are, in addition, robust
with respect to the severe nonlinearities of viscoplastic rheologies and
problem discretizations.
By solvability we refer to the existence of a unique solution of the flow
problem; and convergence in this context means the rate of reduction of the
(nonlinear) residual per iteration.
While a plethora of models exist
that involve plasticity, we focus on a moderately complex but
representative and challenging rheological relation that switches from viscous to
plastic flow at a yield stress. The considered yield stresses may
vary spatially but do not depend on the dynamic pressure and thus not
on the flow solution. The resulting law is sometimes referred to as
von Mises rheology. This law can be generalized in multiple
directions, such as to models that include elasticity, which
naturally introduces a time scale and thus requires time stepping.
The resulting systems exhibit solvability properties similar to
problems without elasticity, but they tend be computationally less challenging
and therefore converge to a solution sufficiently fast using standard nonlinear
solvers \cite{DuretzSoucheBorstEtAl18}.
Another generalization of
viscoplastic models is to problems where the yield stress depends
on the dynamic pressure that is part of the unknown solution, also
called Drucker--Prager models \cite{SpiegelmanMayWilson16,
  PopovSobolev08}. While these models are theoretically and
computationally more complicated, they share several challenges with
von Mises plasticity models, where the yield stress is independent of
the solution.

The available iterative methods for nonlinear Stokes problems can be roughly
divided into two categories: %
Picard-type or
fixed-point methods, and %
Newton-type methods.
It is well documented that %
Picard methods
for nonlinear Stokes equations can be slow to converge and may result in
poorly approximated numerical solutions, even after hundreds of
iterations. Convergence typically degrades with
the ``severity'' of the
nonlinearity, and rheologies involving plasticity are particularly
difficult because of their nonsmoothness.
Thus, in recent years, researchers have been forced either to use
simpler rheology laws or to consider more complex Newton-type or combined
Picard--Newton methods \cite{SpiegelmanMayWilson16,
FratersBangerthThieulotEtAl19, IsaacStadlerGhattas15,
  DuretzSoucheBorstEtAl18, MehlmannRichter17}.
In each iteration of a Newton-type method, a linearized system of equations
is solved numerically to compute a
Newton step, which is subsequently used to improve the current numerical approximation of
the nonlinear system's solution. The derivation and implementation of the linearized
system require %
derivatives/Jacobians,
which can potentially capture the
nonlinear behavior better than Picard methods can, leading to faster
convergence.  For the viscoplastic rheologies we target, however, standard
Newton methods still exhibit slow convergence, which worsens with finer
computational meshes.
Therefore, to advance the body of research on nonlinear methods for
geophysical flows, we propose a novel Newton-type method that requires
solution of Newton linearizations
of the same computational complexity
as standard Newton methods have, and with
similar properties, but that demonstrates
significantly improved convergence properties.
We also
derive
structural properties of non-Newtonian
viscoplastic flow models that we consider relevant for geodynamics. In
particular, we show that underlying two commonly used laws for
viscoplasticity is a minimization problem. This, in turn, provides a
path to solvability, where we show the uniqueness of the solutions by arguing that they
correspond to the unique minimizers of a strictly convex energy functional.

The overall time to solution for a nonlinear Stokes problem comprises the
time for
iterations of the nonlinear solver and for solution of the linearized systems
in each nonlinear step.  Hence, linear iterative methods and preconditioners for
large-scale computational geophysics are important and rich research topics.
The development and implementation of fast and scalable parallel linear solvers
for three-dimensional Stokes problems have been addressed by us
\cite{Rudi18, RudiMalossiIsaacEtAl15, IsaacStadlerGhattas15,
RudiStadlerGhattas17} and others \cite{MayBrownLePourhiet15,
FratersBangerthThieulotEtAl19, MayMoresi08}.
In this paper, however, we focus on nonlinear solvers and refer the reader
to the above literature for numerically solving the linearized systems.
We present the governing equations in Section~\ref{sec:equations} and
show how a viscoplastic flow problem can be derived from energy minimization
principles.
In Section~\ref{sec:optimization}
we discuss what this implies for the uniqueness of solutions and
the need for regularization.
Section~\ref{sec:Newton} presents iterative Newton-type
algorithms, including our novel modified Newton algorithm, in which we
treat the stress as an independent variable in the Newton linearization
but still solve a system only in the standard velocity--pressure
variables. We use two- and three-dimensional benchmark
problems presented in Section~\ref{sec:examples} to numerically study
solution features and to compare the convergence properties
of the proposed solution method with existing approaches. The results
of these comparisons are presented in Section~\ref{sec:results}.

\section{Governing Equations}\label{sec:equations}
We consider nonlinear incompressible Stokes equations on a %
domain $\Omega\subset\mathbb R^{d}$, $d=2,3$, given by
\begin{linenomath}
\begin{subequations}
\label{eq:stokes}
\begin{alignat}{2}
  \label{eq:momentum}
  - \divergence
    \bigl[ \visc(\vel,\press) \, (\gradient\vel + \gradient\vel\transpose) \bigr]
  + \gradient\press &= \rhsvel,
  &&\quad\text{in }\Omega,
  \\
  \label{eq:mass}
  - \divergence\vel &= 0,
  &&\quad\text{in }\Omega,
\end{alignat}
\end{subequations}
\end{linenomath}
where $\vel$ and $\press$ are the velocity and pressure fields,
respectively; the right-hand side volumetric force is denoted by
$\rhsvel$; and the effective viscosity, $\visc$, may depend on the
pressure as well as on the second invariant of the strain rate tensor
$\strainrateinvII^{\vel} \coloneqq (\frac12\,\strainratetensor(\vel) :
\strainratetensor(\vel))^{1/2}$, where ``:'' represents the inner product of
second-order tensors and $\strainratetensor(\vel) \coloneqq \frac12
(\gradient\vel + \gradient\vel\transpose)$ is the strain rate
tensor. We assume that \eqref{eq:stokes} is complemented with
appropriate boundary conditions that, in particular, do not allow for
arbitrary shifts or rotations in the velocity field.
For the discussions in this paper, we further assume that
$\visc=\visc(\vel)$, that is, the viscosity does not depend on
the (dynamic) pressure. Instead, we focus on a
nonsmooth strain-rate dependent constitutive relation that includes a
yield stress, sometimes called the von Mises criterion
\cite{vonMises13, Fullsack95, MoresiSolomatov98,
  SpiegelmanMayWilson16}.  In particular, we consider two
effective viscosity laws for viscoplastic flow,
\begin{linenomath}
\begin{subequations}
\label{eq:viscosity}
\begin{align}
  \label{eq:viscosity-ideal}
  \visc_{i}(\vel) &\coloneqq
    \viscmin + \min\left(
      \frac{\viscstressyield}{2\strainrateinvII^{\vel}},
      \viscref
    \right),\\
  \label{eq:viscosity-composite}
  \visc_{c}(\vel) &\coloneqq
      \frac{\viscstressyield\viscref}{2\strainrateinvII^{\vel}\viscref +
        \viscstressyield}.
  \end{align}
\end{subequations}
\end{linenomath}
These
constitutive relations incorporate plastic yielding with yield
stress $\viscstressyield>0$ and a reference viscosity
$\viscref>0$. Additionally, \eqref{eq:viscosity-ideal} incorporates a
lower bound on the viscosity $\viscmin\ge 0$, which plays a role in
showing that solutions with \eqref{eq:viscosity-ideal} are unique. Both the reference
viscosity and the yield stress can be functions of $\bs x$. In geophysical
applications, $\viscref$ may depend on a spatially varying temperature field,
and the yield stress may depend on both temperature and pressure. We
refer to \eqref{eq:viscosity-ideal} as \emph{ideal} and to
\eqref{eq:viscosity-composite} as \emph{composite} rheology for von
Mises viscoplasticity. The composite rheology law is sometimes
preferred because it avoids the pointwise $\min$-function and thus has
better differentiability properties. Note that $\visc_c$ can also be
written as a scaled harmonic mean, namely, $\visc_c(\vel) =
({2\strainrateinvII^{\vel}}/{\viscstressyield} + 1/\viscref)^{-1}$,
showing the relation between \eqref{eq:viscosity-composite} and
\eqref{eq:viscosity-ideal}.
As part of this paper, we will
also discuss theoretical and practical differences between these two laws.

Next, we illustrate that \eqref{eq:viscosity-ideal} and
\eqref{eq:viscosity-composite} model a viscoplastic fluid with the von
Mises yield criterion; in other words, they model a fluid with a yield
stress. Let us first consider \eqref{eq:viscosity-ideal}, and
introduce the stress tensor $\viscstresstensor\coloneqq 2\min\left(
{\viscstressyield}/{(2\strainrateinvII^{\vel})}, \viscref
\right)\strainratetensor(\vel)$, which allows us to write
\eqref{eq:momentum} in the form
\begin{linenomath}
\begin{align}
  \label{eq:momentum2}
  - \divergence
    \bigl[ 2\viscmin\strainratetensor(\vel) + \viscstresstensor \bigr]
  + \gradient\press &= \rhsvel.
\end{align}
\end{linenomath}
The definition of $\viscstresstensor$ implies that the second
invariant of the stress tensor satisfies
\begin{linenomath}
\begin{equation}
  \viscstressinvII = 2\min\left(
  \frac{\viscstressyield}{2\strainrateinvII^{\vel}}, \viscref\right)
  \strainrateinvII^{\vel} \le \viscstressyield.
  \label{eq:stress-yieldstress-ineq}
\end{equation}
\end{linenomath}
Hence, $\viscstressinvII$ is bounded by the yield stress, and if
$\viscmin=0$, \eqref{eq:viscosity-ideal} models a von Mises rheology with
either constant or spatially varying (e.g., depth-dependent) yield
stress. The parameter $\viscmin\ge 0$ is a regularization parameter such that
$\viscref\gg\viscmin$ and plays an important role in proving uniqueness of a
solution of \eqref{eq:stokes} with \eqref{eq:viscosity-ideal}.

Let us now consider the composite rheology model
\eqref{eq:viscosity-composite}, which is commonly used in geophysical
applications.  Defining the viscous stress tensor as
$\viscstresstensor \coloneqq 2\visc_{c}\strainratetensor(\vel)$, one can easily
show that the second invariant of the corresponding strain rate tensor
satisfies $\viscstressinvII \le \viscstressyield$; that is, it is again
bounded by the yield stress.

The rheologies \eqref{eq:viscosity-ideal} and
\eqref{eq:viscosity-composite} can be generalized in several
directions.  One straightforward generalization is to make
$\viscref$
a function of $\strainrateinvII^{\vel}$, modeling
non-Newtonian behavior in addition to the one caused by the yield
stress.  For instance, polynomial shear thinning is commonly used in
mantle or ice flow \cite{TurcotteSchubert02,Hutter83} models. The
results in this paper can be generalized to these cases, and the only
reason we do discuss models that incorporate
both shear thinning and yielding is clarity of presentation.  In our
mantle flow simulation models \cite{Rudi18,RudiMalossiIsaacEtAl15}, we
combine these two rheological phenomena, generalizing the approach
presented in this paper.  If the yield stress depends on the pressure,
the effective viscosity laws \eqref{eq:viscosity-ideal} and
\eqref{eq:viscosity-composite} are known as Drucker--Prager
rheologies.
\section{Optimization Formulation}\label{sec:optimization}
Here, we characterize the solutions to the incompressible Stokes
equations with either \eqref{eq:viscosity-ideal} or
\eqref{eq:viscosity-composite} as minimizers of an appropriately
chosen energy functional.
The analogous
result for the linear Stokes equations is well known \cite{BrezziFortin91}, namely, that
the linear Stokes solution minimizes a quadratic viscous energy
functional over a space of sufficiently regular divergence-free
vector functions.
For simplicity, here we  assume that the
Stokes systems are combined with homogeneous Dirichlet boundary
conditions on the entire boundary of the domain $\partial\Omega$. Other
boundary conditions, such as free-slip conditions, can be used as long as they
eliminate the null
space of the Stokes operator, that is, exclude translations or rotations
from the solution space \cite{ElmanSilvesterWathen14a}.

\begin{figure}\centering
	\includegraphics[width=0.4\columnwidth]{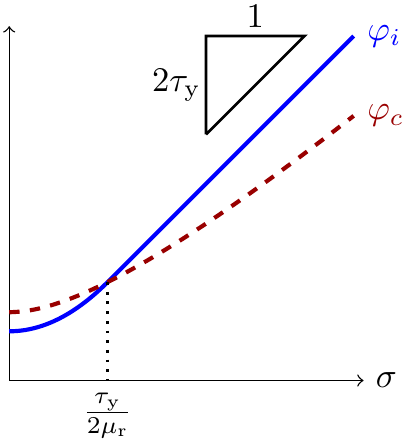}
  \caption{Sketch of functions $\varphi_i$ and $\varphi_c$ for ideal
    and composite viscoplasticity, respectively, as defined in \eqref{eq:Phi}. The
    function $\varphi_i$ consists of a quadratic and a linear piece,
    which are joined at
    $\sigma=\viscstressyield/(2\viscref)$. Both functions approach a slope
    of $2\viscstressyield$ as $\sigma\to\infty$.
  }\label{fig:Phi}
\end{figure}

To state the energy functional for the Stokes
problem \eqref{eq:stokes} and \eqref{eq:viscosity}, we first define
continuously differentiable functions $\varphi:\mathbb R_{\ge 0}\to
\mathbb R$ representing ideal and composite viscoplasticity as follows:
\begin{linenomath}
\begin{subequations}\label{eq:Phi}
\begin{align}\label{eq:Phi-ideal}
\varphi_i(\sigma) &\coloneqq
\begin{cases}
	2\viscref\sigma^2 + \dfrac{\viscstressyield^2}{2\viscref} & \text{ if }
	2\viscref\sigma\le \viscstressyield,\\
    2\viscstressyield \sigma & \text{ else,}
  \end{cases}\\
\label{eq:Phi-composite}
\varphi_c(\sigma) &\coloneqq
2\viscstressyield\sigma - \dfrac{\viscstressyield^2}{\viscref}
\log(\viscstressyield+2\viscref\sigma),
\end{align}
\end{subequations}
\end{linenomath}
where as before we assume
$\viscref,\viscstressyield>0$. Graphs of these functions are given
in Figure~\ref{fig:Phi} for illustration. We consider the energy minimization
\begin{linenomath}
\begin{equation}\label{eq:minJ}
  \min_{\nabla\cdot{\bs u} = 0} \Phi(\bs u),\quad
  \Phi({\bs u}) \coloneqq \int_\Omega \varphi(\strainrateinvII^{\vel})\,d\bs
  x + 2\viscmin\int_\Omega (\strainrateinvII^{\vel})^2\,d\bs
  x - \int_\Omega \rhsvel\cdot {\bs u}\,d\bs x,
\end{equation}
\end{linenomath}
where $\varphi$ can be either $\varphi_i$ or $\varphi_c$. In the
former case, in which we consider the function corresponding to the ideal
law \eqref{eq:Phi-ideal}, we
assume $\viscmin\ge 0$.  In the case in which we consider the composite
rheology function $\varphi=\varphi_c$, we always assume $\viscmin=0$
such that the quadratic term in $\Phi(\cdot)$
vanishes. The
minimization in \eqref{eq:minJ} is over all divergence-free functions
${\bs u}\in H_0^1(\Omega)^d$.  Here, the function space
$H_0^1(\Omega)^d$ contains square-integrable functions that have
square-integrable derivatives and that (for simplicity) satisfy
homogeneous Dirichlet boundary conditions on $\partial\Omega$.

Next, we compute the first-order necessary conditions
for \eqref{eq:minJ}, that is, the conditions a minimizer of
\eqref{eq:minJ} must satisfy.  For that purpose, we use an arbitrary
divergence-free
velocity $\bs v\in H_0^1(\Omega)^d$ that satisfies the same
homogeneous boundary condition as ${\bs u}$. Recall that the
variation of $\strainrateinvII^{\bs u}$ with respect to $\bs u$ in a direction
$\bs v$ is given by
\begin{linenomath}
\begin{equation}\label{eq:derivII}
  [\delta_{{\bs u}}\strainrateinvII^{\bs u}](\bs v) =
    \biggl[\delta_{{\bs u}}
          \sqrt{\tfrac12\,\strainratetensor(\vel) : \strainratetensor(\vel)}\biggr](\bs v) =
    \frac{\strainratetensor(\bs u) : \strainratetensor(\bs v)}{2\strainrateinvII^{\bs u}}.
\end{equation}
\end{linenomath}
Here and in the following, $\delta_{\vel}$ denotes the variation with
respect to $\vel$.
For the variation of $\Phi$ at $\bs u$ in direction $\bs v$, using the chain
rule and \eqref{eq:derivII}, we obtain
\begin{linenomath}
\begin{equation}\label{eq:deriv-Phi}
  [\delta_{{\bs u}} \Phi({\bs u})](\bs v) = \int_\Omega
  \varphi'(\strainrateinvII^{\bs u}) \frac{\strainratetensor(\bs{u}) :
    \strainratetensor(\bs{v})}{2\strainrateinvII^{\bs u}} \,d\bs
  x + \int_{\Omega} 2\viscmin\left(\strainratetensor(\bs{u}) :
  	\strainratetensor(\bs{v})\right) \,d\bs
  x - \int_\Omega \rhsvel\cdot \bs v\,d\bs x.
\end{equation}
\end{linenomath}
Taking the derivatives of \eqref{eq:Phi-ideal} and \eqref{eq:Phi-composite},
we have
\begin{linenomath}
\begin{align}
  \varphi_i'(\sigma) &=
  \left.\begin{cases}
    4\viscref \sigma  & \text{ if }
    2\viscref\sigma< \viscstressyield\\
    2\viscstressyield & \text{ if } 2\viscref\sigma> \viscstressyield
  \end{cases}\right\} = 2\min(2\viscref\sigma,\viscstressyield),\\
  \varphi_c'(\sigma) &= \frac{4\viscref\viscstressyield\sigma}{2\viscref\sigma+\viscstressyield},
\end{align}
\end{linenomath}
we obtain for the integrands in \eqref{eq:deriv-Phi}
\begin{linenomath}
\begin{align*}
\varphi_i'(\strainrateinvII^{\bs u}) \frac{\strainratetensor(\bs u)
  : \strainratetensor(\bs{v})}{2\strainrateinvII^{\bs{u}}}
&=  \left[ \min\left(2\viscref,\dfrac{\viscstressyield}{\strainrateinvII^{\bs u}}\right)\right]
  \strainratetensor(\bs{u}) : \strainratetensor(\bs{v}),\\
\varphi_c'(\strainrateinvII^{\bs u}) \dfrac{\strainratetensor(\bs{u})
	: \strainratetensor(\bs{v})}{2\strainrateinvII^{\bs u}}
&=  \left[\dfrac{2\viscref\viscstressyield}{2\viscref\strainrateinvII^{\bs u}+\viscstressyield}\right]
\strainratetensor(\bs{u}) : \strainratetensor(\bs{v}).
\end{align*}
\end{linenomath}
Hence, the necessary optimality condition, $[\delta_{{\bs u}} \Phi({\bs
  u})](\bs v)=0$ for all admissible $\bs v$, is the variational form of the
the strong form of the Stokes problem \eqref{eq:stokes}, \eqref{eq:viscosity}. The strong
form can be recovered using integration by parts and the fact that
$\bs v$ is arbitrary and by introducing the pressure as a Lagrange multiplier
for the divergence-free condition \cite{BrezziFortin91}.

Note that the optimization functional \eqref{eq:minJ} is convex, which
makes it useful for arguing existence and uniqueness of solutions. In
particular if $\viscmin>0$, the Stokes problem \eqref{eq:stokes} and
\eqref{eq:viscosity}, combined with appropriate boundary conditions, has
a unique solution since the energy functional $\Phi$ in \eqref{eq:minJ}
is strictly convex.
While for $\viscmin=0$ existence of solutions can
be argued, they are not necessarily unique. Additionally, since the
Stokes solution corresponds to a minimum of $\Phi$, this energy
functional can also be used in a line search procedure to guarantee
descent (and thus convergence) of iterative methods for nonlinear systems, such as
a Newton or Picard method.

\section{Standard and Stress--Velocity Newton Methods}\label{sec:Newton}
Common methods for solving nonlinear systems of equations of the form
\eqref{eq:stokes} are the Picard fixed-point method and Newton's method, or
combinations thereof \cite{MoresiSolomatov98, Rudi18, FratersBangerthThieulotEtAl19,
  MayBrownLePourhiet15, DuretzSoucheBorstEtAl18, aspect,
  SpiegelmanMayWilson16, MehlmannRichter17}. As discussed
in the preceding section, these iterative methods can be seen either as
solving the nonlinear Stokes equations \eqref{eq:stokes} or as
minimizing the energy functional \eqref{eq:minJ}. Both iterative
approaches are summarized next, where our derivations initially use
an unspecified (convex) function $\varphi$ as in \eqref{eq:minJ}.
At the end of this section, we list the standard and the stress--velocity
Newton's methods in Algorithm \ref{alg:newton} as well as the specific
linearized systems for ideal viscoplasticity $\varphi=\varphi_i$ in Table
\ref{tab:ideal} and for composite viscoplasticity $\varphi=\varphi_c$ in Table
\ref{tab:composite}.

\subsection{Picard and standard Newton methods}
The Picard fixed-point iteration for
\eqref{eq:stokes} and \eqref{eq:viscosity} starts with an initial guess
$(\vel_0,p_0)$ for velocity and pressure and, given an iterate
$(\vel_{k-1},p_{k-1})$, computes the next iterate $(\vel_k,p_k)$ for
$k=1,2,\ldots$ as solution to
\begin{linenomath}
\begin{subequations}\label{eq:picard}
\begin{alignat}{1}%
  \label{eq:picard:momentum}
  - \nabla\cdot
  \left[\left(2\viscmin + \frac{\varphi'(\strainrateinvII^{\bs u_{k-1}})}{2\strainrateinvII^{\bs u_{k-1}}}\right)
   \strainratetensor(\bs{u_k} ) \right]
  + \nabla p_k &= \rhsvel,
  \\
  \label{eq:picard:mass}
  -\nabla\cdot\vel_k &= 0,
\end{alignat}
\end{subequations}
\end{linenomath}
where
$\strainrateinvII^{\bs u_{k-1}}$ is computed from
the ($k-1$)st iterate.
The new iterate depends only on the previous velocity and
not on the pressure since the pressure variable enters linearly and
acts as the Lagrange multiplier for enforcing the divergence-free
condition.  The Picard fixed-point method requires, in each iteration,
the solution of one linear Stokes equation with spatially varying but
isotropic viscosity. It typically converges from any initialization
$(\vel_0,p_0)$,  but often exhibits slow convergence,
in particular in the
presence of significantly nonlinear or nonsmooth viscosity laws such as \eqref{eq:viscosity}.

An alternative to the Picard iteration is Newton's method, which,
starting from an initial velocity--pressure pair $(\vel_0,p_0)$,
computes iterates by first solving the linear Newton system for the
Newton update $(\tilde\vel,\tilde p)$:
\begin{linenomath}
\begin{subequations}\label{eq:primal:newton}
  \begin{alignat}{2}
\label{eq:primal:momentum}
- \nabla\cdot
\left[2\viscmin\strainratetensor(\tilde{\vel})+\frac{\varphi'(\strainrateinvII^{\bs u_{k-1}})}{2\strainrateinvII^{\bs u_{k-1}}}
\bigg(\mathbb I - \psi(\strainrateinvII^{\bs u_{k-1}}) \dfrac{\strainratetensor(\vel_{k-1})\otimes \strainratetensor(\vel_{k-1})}{2(\strainrateinvII^{\bs u_{k-1}})^2} \bigg)
\strainratetensor(\tilde{\vel})\right]
+ \nabla \tilde p &= -\bs r^{\vel}_{k-1},
\\
\label{eq:primal:mass}
-\nabla\cdot\tilde{\bs u} &= -r^{p}_{k-1},
\end{alignat}
where we use the notation
\begin{equation}\label{eq:psi}
\psi(\strainrateinvII^{\bs u_{k-1}}) \coloneqq \dfrac{\varphi'(\strainrateinvII^{\bs u_{k-1}})-\strainrateinvII^{\bs u_{k-1}}\varphi''(\strainrateinvII^{\bs u_{k-1}}) }{\varphi'(\strainrateinvII^{\bs u_{k-1}})}.
\end{equation}
\end{subequations}
\end{linenomath}
This is followed by the Newton update with step length $\alpha\le 1$:
\begin{linenomath}
\begin{equation}\label{eq:primal:update}
  \vel_k = \vel_{k-1} + \alpha\tilde\vel, \quad
  p_k = p_{k-1} + \alpha\tilde p,
\end{equation}
\end{linenomath}
which for $\alpha=1$ amounts to a so-called full Newton step.
In \eqref{eq:primal:momentum}, $\mathbb I$
is the 4th-order identity tensor, and ``$\otimes$'' is the outer product
between two (second-order) strain rate tensors.
In the derivation of \eqref{eq:primal:momentum}, we have used the identity
$(\bs a:\bs b)\bs c = (\bs b\otimes \bs c)\bs a $ for second-order tensors
$\bs a, \bs b, \bs c$.
The right-hand sides of the Newton system, $\bs r^{\vel}_{k-1}:=\bs r^{\vel}(\vel_{k-1},p_{k-1})$ and
$r^{p}_{k-1}:=r^p(\vel_{k-1})$, are the nonlinear residuals of the momentum
\eqref{eq:momentum} and mass \eqref{eq:mass} equations, respectively.  Each Newton iteration requires
the solution of a linear Stokes equation with a 4th-order anisotropic
viscosity tensor, which is
symmetric and positive semidefinite. The latter can be seen
directly for $\varphi=\varphi_i$ and $\varphi=\varphi_c$ (see Tables
\ref{tab:ideal} and \ref{tab:composite}) but also follows from the
convexity properties of the objective \eqref{eq:minJ}. If
$\varphi=\varphi_c$, the 4th-order tensor in \eqref{eq:primal:newton}
is positive definite, and thus this Newton linearization has a unique
solution. If $\varphi=\varphi_i$, then this tensor is only
semidefinite; but then assuming that $\viscmin>0$ guarantees that the
Newton system has a unique solution.

Note that the Picard method is a special case of Newton's method. It
arises from the Newton system \eqref{eq:primal:newton} by neglecting the
higher-order derivatives, that is, setting
$\psi(\strainrateinvII^{\bs u_{k-1}})= 0$,
and by substituting $(\tilde \vel,\tilde p)$ in
\eqref{eq:primal:newton} using
\eqref{eq:primal:update} with $\alpha=1$:
\begin{linenomath}
\begin{align*}
  - \nabla\cdot
  \left[\left(2\viscmin + \frac{\varphi'(\strainrateinvII^{\bs u_{k-1}})}{2\strainrateinvII^{\bs u_{k-1}}}\right)
   \left(\strainratetensor(\bs{u_k}) - \strainratetensor(\bs{u_{k-1}})\right) \right]
  + \nabla \left(p_k - p_{k-1}\right)
  &= -\bs r^{\vel}_{k-1},
  \\
  -\nabla\cdot \left(\vel_k - \vel_{k-1}\right) &= -r^{p}_{k-1}.
\end{align*}
\end{linenomath}
Recognizing that the negative residuals
are the momentum and mass equations evaluated at the $(k-1)$st
iterate, we obtain the Picard system \eqref{eq:picard}.
As
documented in the literature,
Picard iterations exhibit in general better robustness compared with Newton's
method.  To take
advantage of Newton's faster convergence locally near the solution,
a number of Picard iterations are sometimes run before
switching to Newton's method
\cite{SpiegelmanMayWilson16,MayBrownLePourhiet15}.  There is, however,
no systematic way to determine at which iteration to switch from
Picard to Newton, and one is left with heuristics from trial and error.
The next section presents an alternative Newton linearization, which aims to
combine the robustness of Picard's method with the local fast convergence of
Newton's method.  The method requires the
solution of a modified system that has a structure similar to the standard
approach \eqref{eq:primal:newton}.

\subsection{A stress--velocity Newton method}\label{sec:pdNewton}
Here, we present an alternative formulation of the nonlinear Stokes
problem \eqref{eq:stokes} and \eqref{eq:viscosity}, which leads to a
modified Newton algorithm. This formulation is motivated by the
optimization formulation \eqref{eq:minJ} and  primal-dual
Newton methods in rather different application areas such as
total-variation image denoising \cite{ChanGolubMulet99,
  HintermullerStadler06}, linear elasticity problems with friction and
plasticity \cite{HueberStadlerWohlmuth08, HagerWohlmuth09} or Bingham
fluids \cite{DelosreyesGonzalez09}.  We introduce an
independent variable for the viscous stress tensor as in
\eqref{eq:momentum2} and use that the following two conditions are
equivalent:
\begin{linenomath}
\begin{equation}\label{eq:NCPs}
\viscstresstensor - \frac{\varphi'(\strainrateinvII^{\bs u})}{2\strainrateinvII^{\bs u}}
\strainratetensor(\bs{u}) = 0 \quad \Leftrightarrow \quad
\frac{2\strainrateinvII^{\bs u}}{\varphi'(\strainrateinvII^{\bs u})}\viscstresstensor - \strainratetensor(\bs{u}) = 0.
\end{equation}
\end{linenomath}
Note that this equivalency holds as
${\varphi'(\strainrateinvII^{\bs u})}/{2\strainrateinvII^{\bs u}}  =
2\visc(\vel)>0$ provided that $\strainrateinvII^{\bs u}$
is finite.
The resulting stress--velocity formulation for the unknowns $(\bs u, p, \bs
\tau)$ is derived by substitution of \eqref{eq:NCPs} into \eqref{eq:stokes} and
adding \eqref{eq:NCPs} as a new equation:
\begin{linenomath}
\begin{subequations}\label{eq:pd:stokes}
\begin{alignat}{2}
\label{eq:pd:momentum}
- \divergence\left(2\viscmin\strainratetensor(\bs{u})+
\viscstresstensor\right)
+ \gradient\press &= \rhsvel,
&&\quad\text{in }\Omega,
\\
\label{eq:pd:mass}
-\nabla\cdot{\bs u} &= 0,
&&\quad\text{in }\Omega,
\\
\label{eq:pd:tau}
\frac{2\strainrateinvII^{\bs u}}{\varphi'(\strainrateinvII^{\bs u})}\viscstresstensor - \strainratetensor(\bs{u}) &= 0,
&& \quad\text{in }\Omega.
\end{alignat}
\end{subequations}
\end{linenomath}
Before computing the linearization of the full nonlinear system \eqref{eq:pd:stokes}, note
that only \eqref{eq:pd:tau} is nonlinear, for which we introduce the
notations
\begin{linenomath}
\begin{equation}\label{eq:residual-tau}
\bs{r}^{\bs \tau}(\bs u,\bs \tau): = \frac{2\strainrateinvII^{\bs \vel}}{\varphi'(\strainrateinvII^{\bs \vel})}\viscstresstensor - \strainratetensor(\bs{\vel}).
\end{equation}
\end{linenomath}
To compute directional derivatives of $\bs r^{\bs \tau}$ with respect to
$(\vel,\viscstresstensor)$, we denote the corresponding variations by
$(\tilde{\vel},\tilde{\viscstresstensor})$ and use \eqref{eq:derivII} to obtain
\begin{linenomath}
\begin{equation}\label{eq:residual-tau-deriv}
  [\delta \bs{r}^{\bs \tau}(\vel,\bs \tau)](\tilde{\vel},\tilde{\viscstresstensor}) =
  \dfrac{2\varphi'(\strainrateinvII^{\vel})-2\strainrateinvII^{\vel}\varphi''(\strainrateinvII^{\vel}) }
        {\varphi'(\strainrateinvII^{\vel})^2}
  \dfrac{\strainratetensor(\vel) : \strainratetensor(\tilde{\vel})}{2\strainrateinvII^{\vel}}
  \viscstresstensor +
  \dfrac{2\strainrateinvII^{\vel}}{\varphi'(\strainrateinvII^{\vel})} \tilde{\viscstresstensor} -
  \strainratetensor(\tilde{\vel}).
\end{equation}
\end{linenomath}
Therefore, the Newton step for the full nonlinear system \eqref{eq:pd:stokes} is
\begin{linenomath}
\begin{subequations}\label{eq:pd:stokes-deriv}
\begin{alignat}{2}
\label{eq:pd:momentum-deriv}
- \divergence\left(2\viscmin\strainratetensor(\tilde{\vel})+
\tilde{\viscstresstensor}\right)
+ \gradient\tilde{\press} &= \rhsvel +
\divergence\left(2\viscmin\strainratetensor(\bs{u})+ \viscstresstensor\right)
- \gradient\press
\\
\label{eq:pd:mass-deriv}
-\nabla\cdot{\tilde{\vel}} &= \nabla\cdot{\bs u}
\\
\label{eq:pd:tau-deriv}
[\delta \bs{r}^{\bs \tau}(\bs u,\bs \tau)](\tilde{\bs u},\tilde{\bs \tau}) &= -\bs{r}^{\bs \tau}(\bs u,\bs \tau),
\end{alignat}
\end{subequations}
\end{linenomath}
Next, we eliminate \eqref{eq:pd:tau-deriv} from this system
by isolating
\begin{linenomath}
\begin{equation}\label{eq:tau-var}
\tilde{\viscstresstensor}
=
\left(\dfrac{\varphi'(\strainrateinvII^{\vel})}{2\strainrateinvII^{\vel}}-
\dfrac{\psi(\strainrateinvII^{\vel}) }
{\strainrateinvII^{\vel}}
\dfrac{\strainratetensor(\vel)\otimes \viscstresstensor}{2\strainrateinvII^{\vel}}
\right)\strainratetensor(\tilde{\vel}) -\dfrac{\varphi'(\strainrateinvII^{\vel})}{2\strainrateinvII^{\vel}}\bs{r}^{\bs \tau}(\bs u,\bs \tau),
\end{equation}
\end{linenomath}
where ``$\otimes$'' again denotes the outer product between two second-order
tensors and $\psi(\cdot)$ is defined as in \eqref{eq:psi}.
Substitution of $\tilde{\viscstresstensor}$ in \eqref{eq:tau-var} into
\eqref{eq:pd:momentum-deriv} eliminates \eqref{eq:pd:tau-deriv} and yields the
following reduced Newton system (at Newton iteration $k$):
\begin{linenomath}
\begin{subequations}\label{eq:pd1:stokes}
	\begin{alignat}{2}
	\label{eq:pd1:momentum}
	- \nabla\cdot
	\left[2\viscmin\strainratetensor(\tilde{\vel}) + \frac{\varphi(\strainrateinvII^{\bs u_{k-1}})}{2\strainrateinvII^{\bs u_{k-1}}}
	\bigg(\mathbb I -
	\dfrac{\sqrt{2} \psi(\strainrateinvII^{\bs u_{k-1}})}{\varphi'(\strainrateinvII^{\bs u_{k-1}})}
	\dfrac{\strainratetensor(\vel_{k-1})\otimes \viscstresstensor_{k-1}}{\sqrt{2}\strainrateinvII^{\vel_{k-1}}} \bigg)
	\strainratetensor(\tilde{\vel})\right]\!
	+ \!\nabla \tilde p &= -\bs r^{\vel}_{k-1},
	\\
	\label{eq:pd1:mass}
	-\nabla\cdot\tilde{\bs u} &= -r^{p}_{k-1},
\end{alignat}
\end{subequations}
\end{linenomath}
where the right hand sides are, as before, the negative residuals of
\eqref{eq:momentum} and \eqref{eq:mass}.
Thus, to compute the full Newton step
$(\tilde{\bs u},\tilde{p},\tilde{\bs \tau})$, only the reduced system
\eqref{eq:pd1:stokes} needs to be solved for $(\tilde{\bs u},\tilde{p})$
numerically, which amounts to the same computational complexity as computing
a step of the standard Newton linearization \eqref{eq:primal:newton}.  This
is a key property of our stress--velocity Newton formulation because the
overall time to (nonlinear) solution of the two Newton linearizations can
simply be compared by counting the number of Newton iterations.  After the
numerical solution of \eqref{eq:pd1:stokes} for $(\tilde{\bs u},\tilde{p})$ is
available, one can,
evaluate $\tilde{\viscstresstensor}$ by replacing $\vel$ by $\vel_{k-1}$ and
$\viscstresstensor$ by $\viscstresstensor_{k-1}$ in \eqref{eq:tau-var}.
The computational cost of this evaluation is negligible compared with the cost
of one linear solve.

Next, we ask whether the modified Newton step \eqref{eq:pd1:stokes}
always has unique solution. We will find that this can
only be guaranteed after a small modification that is motivated by the
optimization formulation and by comparison with the standard Newton
method. However, this modification is crucial for robust convergence
in numerical practice.
Since we introduced the
additional stress variable $\viscstresstensor$ during the
linearization, we have lost the direct connection to an optimization
problem and thus cannot use convexity arguments for an appropriate
objective function.
Comparing the modified Newton step
\eqref{eq:pd1:stokes} with the standard Newton step
\eqref{eq:primal:newton}, we find that the only difference is the
replacement of
\begin{linenomath}
\begin{equation*}
  \frac{\psi(\strainrateinvII^{\bs u_{k-1}})}
       {\sqrt{2}\strainrateinvII^{\bs u_{k-1}}}
       \strainratetensor
  \text{ in } \eqref{eq:primal:momentum}
  \quad\text{by}\quad
  \frac{\sqrt{2}\psi(\strainrateinvII^{\bs u_{k-1}})}
       {\varphi'(\strainrateinvII^{\bs u_{k-1}})}
       \viscstresstensor
  \text{ in } \eqref{eq:pd1:momentum}.
\end{equation*}
\end{linenomath}
Since the second invariant of the left tensor is bounded by 1,
we must ensure that the second invariant of the corresponding
tensor in the stress--velocity formulation is also bounded by 1. This
is the case if $\viscstressinvII \le \viscstressyield$,
namely, that the stress tensor satisfies the yield stress bound. Note that
this inequality is automatically satisfied when $\bs r^{\bs \tau}(\vel,
\viscstresstensor)=\boldsymbol 0$, according to
\eqref{eq:stress-yieldstress-ineq}.  Unfortunately, this
does not hold generally. Thus, we ensure this by modifying the
term involving the stress in \eqref{eq:pd1:momentum} by scaling it
and using a symmetrized form of the 4th-order tensor. This
rescaling is of the form
\begin{linenomath}
\begin{equation}\label{eq:tau_projection}
\bs \tau \leftarrow \frac{ \bs \tau}{\max(1,\viscstressinvII/\viscstressyield)}.
\end{equation}
\end{linenomath}
Combined with the symmetrization, this results in replacing
\begin{linenomath}
\begin{equation}\label{eq:tau-modification}
\dfrac{\strainratetensor(\vel_{k-1})\otimes
  \viscstresstensor_{k-1}}{\sqrt{2}\strainrateinvII^{\vel_{k-1}}}
\quad \text{ by } \quad
\dfrac{\left(\strainratetensor(\vel_{k-1})\otimes
 \viscstresstensor_{k-1}\right)_{\textit{sym}}
 }{\sqrt{2}\strainrateinvII^{\vel_{k-1}}\max(1,\viscstressinvII/\viscstressyield)},
\end{equation}
\end{linenomath}
where $\left(\bs a\otimes \bs b\right)_{\textit{sym}} \coloneqq \frac 1 2 \left(\bs a\otimes \bs b + \bs b\otimes \bs a\right)$ is the symmetrization of the 4th-order tensor $\bs a\otimes \bs b$.

This modification to the stress--velocity Newton problem makes the
algorithm a quasi-Newton method, i.e., the exact derivative matrix is
replaced with an approximation. Upon convergence of the
algorithm, however, the modification \eqref{eq:tau-modification} becomes
smaller since in the limit, $\bs r^{\bs \tau}(\vel, \viscstresstensor)\rightarrow0$,
\eqref{eq:NCPs} must hold. This implies
that, in the limit, the second invariant of the stress tensor is
bounded by the yield stress and that the stress tensor is a multiple
of the strain tensor and thus the influence of the symmetrization in
\eqref{eq:tau-modification} decreases. Since upon convergence the
modifications vanish, we can expect to recover fast Newton-type
convergence close to the solution. Note that similar modifications for
primal-dual Newton methods, applied to rather different problems, are
also discussed in \citeA{ChanGolubMulet99, HintermullerStadler06,
  HueberStadlerWohlmuth08, HagerWohlmuth09, DelosreyesGonzalez09}.

Moreover, note that when $\viscstresstensor=\boldsymbol 0$, \eqref{eq:pd1:stokes}
simply takes the form of the Picard fixed-point method. Since Picard
is known to be stable and also to converge rapidly at the beginning (i.e., far from the solution),
$\viscstresstensor=\boldsymbol 0$ is a good initialization. After initial
iterations, the stress--velocity Newton method becomes increasingly similar and
eventually converges to the
standard Newton method. This can be seen by replacing
$\viscstresstensor_{k-1}$ in \eqref{eq:pd1:momentum} by
$\varphi'(\strainrateinvII^{\vel_{k-1}})/(2\strainrateinvII^{\vel_{k-1}})\strainratetensor(\vel_{k-1})$
following the definition of the stress tensor \eqref{eq:NCPs}.
This identity is satisfied upon convergence of
$\bs r^{\bs \tau}(\vel,\viscstresstensor)\to \boldsymbol 0$.
\subsection{Summary of methods for ideal and composite viscoplasticity}
We summarize the methods developed above and specify the
algorithmic steps  for the ideal and composite rheological laws. For
this purpose, we first collect derivatives of the functions
$\varphi_i$ and $\varphi_c$:
\vskip 1.5ex
\begin{tabular}{ll}
$\varphi_i(\sigma) =
\begin{cases}
2\viscref \sigma^2 + \dfrac{\viscstressyield^2}{2\viscref} & \text{ if }
2\viscref\sigma\le \viscstressyield,\\
2\viscstressyield \sigma & \text{ else}
\end{cases}$ &
$\varphi_c(\sigma) = 2\viscstressyield\sigma - \dfrac{\viscstressyield^2}{\viscref} \log(\viscstressyield+2\viscref\sigma)$\\

\\

$\varphi_i'(\sigma) = \begin{cases}
4\viscref \sigma & \text{ if }
2\viscref\sigma < \viscstressyield,\\
2\viscstressyield& \text{ if } 2\viscref\sigma > \viscstressyield
\end{cases}$ &
$\varphi_c'(\sigma) = \dfrac{4\viscref\viscstressyield\sigma}{2\viscref\sigma+\viscstressyield}$\\

\\

$\varphi_i''(\sigma) = \begin{cases}
4\viscref & \text{ if }
2\viscref\sigma < \viscstressyield,\\
0& \text{ if } 2\viscref\sigma > \viscstressyield
\end{cases}$ &
$\varphi_c''(\sigma) =
\dfrac{4\viscref\viscstressyield^2}{\left(2\viscref\sigma+\viscstressyield\right)^2}$\\
\end{tabular}
\vskip 1.5ex

The Picard, standard Newton, and stress--velocity Newton
problems are summarized in Table \ref{tab:ideal} for ideal
viscoplasticity and in Table \ref{tab:composite} for composite
viscoplasticity.  In these tables, the Picard method is written
directly for the new iterate $(\vel_k,p_k)$, whereas the Newton
methods are written for the updates $(\tilde\vel,\tilde p)$; then
the
iterate $(\vel_{k},p_k)$ is computed as $(\vel_{k},p_k) =
(\vel_{k-1},p_{k-1}) + \alpha (\tilde\vel,\tilde p)$ with a step
length $\alpha$ that is usually equal to or less than 1.
After computing the stress--velocity Newton
update, the corresponding stress update is found by using
\eqref{eq:tau-var}. In the stress--velocity Newton updates, the term
involving $\viscstresstensor_{k-1}$ is modified as shown in
\eqref{eq:tau-modification} to ensure solvability of the
stress--velocity Newton step. For a complete listing of standard and
stress--velocity Newton's methods see Algorithm \ref{alg:newton}, where the
labels \textbf{(N)}, \textbf{(SVN)}, and \textbf{(SU)} refer to either of the
equations in Table \ref{tab:ideal} or \ref{tab:composite} depending whether
the ideal or the composite viscoplastic law is used.

\begin{algorithm}
\caption{Standard and stress--velocity Newton's methods}
\label{alg:newton}
\begin{minipage}[t]{0.49\columnwidth}
  \vspace{0.5ex}
  \textbf{Standard Newton's method:}
  \begin{algorithmic}[1]
    \State Choose initial guess $(\vel_0,p_0)$
    \For{$k = 1, 2, \ldots$ until convergence}
      \State Evaluate residual $(\bs r^{\vel}_{k-1},r^{p}_{k-1})$
      \State Solve system \textbf{(N)} for %
             $(\tilde\vel,\tilde p)$
      \State Find step length $\alpha$ via line search
      \State Set %
             $\vel_k = \vel_{k-1}+\alpha\tilde\vel$,
             $p_k = p_{k-1}+\alpha\tilde p$
    \EndFor
  \end{algorithmic}
  \vspace{1.0ex}
\end{minipage}
\hfill
\begin{minipage}[t]{0.49\columnwidth}
  \vspace{0.5ex}
  \textbf{Stress--velocity Newton's method:}
  \begin{algorithmic}[1]
    \State Choose initial guess $(\vel_0,p_0)$, set $\viscstresstensor_0=\boldsymbol 0$
    \For{$k = 1, 2, \ldots$ until convergence}
      \State Evaluate residual $(\bs r^{\vel}_{k-1},r^{p}_{k-1})$
      \State Solve system \textbf{(SVN)} for %
             $(\tilde\vel,\tilde p)$
      \State Evaluate stress update $\tilde{\viscstresstensor}$ as in \textbf{(SU)}
      \State Find step length $\alpha$ via line search
      \State Set %
             $\vel_k = \vel_{k-1}+\alpha\tilde\vel$,
             $p_k = p_{k-1}+\alpha\tilde p$,
      \Statex\hspace{3em}
             $\viscstresstensor_k = \viscstresstensor_{k-1}+
                                    \alpha\tilde{\viscstresstensor}$
    \EndFor
  \end{algorithmic}
  \vspace{1.0ex}
\end{minipage}
\end{algorithm}

\begin{table}
\centering
\caption{Comparison of iterative schemes for \emph{ideal viscoplastic law}
  \eqref{eq:viscosity-ideal}. $\mathcal X$ is the indicator function of the set $\{2\strainrateinvII^{\vel_{k-1}}\viscref > \viscstressyield\}$. %
}\label{tab:ideal}
\begin{tabular}{l}
  \toprule
Picard solve \textbf{(P)}:\\[1ex]
{$\begin{aligned}
- \nabla\cdot
\left[\left(2\viscmin + \min\Big(2\viscref,
\frac{\viscstressyield}{\strainrateinvII^{\bs u_{k-1}}}\Big)\right)
\strainratetensor(\bs u_k) \right] + \nabla p_k &= \rhsvel\\
-\nabla\cdot\vel_k &= 0\\
\end{aligned}$}\\
\midrule
Standard Newton solve \textbf{(N)}:\\[1ex]
{$\begin{aligned}
	- \nabla\cdot
	\left[\left(2\viscmin + \min\Big(2\viscref,
	\frac{\viscstressyield}{\strainrateinvII^{\vel_{k-1}}}\Big)
	\bigg(\mathbb I - \mathcal X
	\dfrac{\strainratetensor(\vel_{k-1})\otimes \strainratetensor(\vel_{k-1})}
	{2(\strainrateinvII^{\vel_{k-1}})^2} \bigg)\right)
	\strainratetensor(\tilde{\vel})\right]
	+ \nabla \tilde p &= -\bs r^{\vel}_{k-1}
	\\
	-\nabla\cdot\tilde{\bs u} &= -r^{p}_{k-1}
	\end{aligned}$}\\[1ex]
\midrule
Stress--velocity Newton solve \textbf{(SVN)}:\\[1ex]
{$
\begin{aligned}
  - \nabla\cdot
   \left[\left(2\viscmin+\min\left(2\viscref,
   \frac{\viscstressyield}{\strainrateinvII^{\vel_{k-1}}}
   \right)\left(\mathbb I  -  \mathcal X
    \dfrac{\left(\strainratetensor(\vel_{k-1})\otimes \viscstresstensor_{k-1}\right)_{\textit{sym}}}
    {2\strainrateinvII^{\vel_{k-1}}\max(\viscstressyield,\viscstressinvII)} \right)\right)
    \strainratetensor(\tilde{\vel})\right]
  + \nabla \tilde p &= -\bs r^{\vel}_{k-1}
  \\
  -\nabla\cdot\tilde{\bs u} &=-r^{p}_{k-1}\\
\end{aligned}$}\\
Stress update \textbf{(SU)}:\\
{$
  \tilde{\viscstresstensor} =-\viscstresstensor_{k-1} + \min\left(2\viscref,
  \frac{\viscstressyield}{\strainrateinvII^{\vel_{k-1}}}
  \right)\left(\strainratetensor(\tilde{\vel}) + \strainratetensor(\vel_{k-1})\right)-\mathcal X
  \dfrac{\left(\strainratetensor(\vel_{k-1})\otimes \viscstresstensor_{k-1}\right)_{\textit{sym}}}{2\left(\strainrateinvII^{\vel_{k-1}}\right)^2\max(1,\viscstressinvII/\viscstressyield)}
  \strainratetensor(\tilde{\vel})
$}\\
\bottomrule
\end{tabular}
\end{table}

\begin{table}
\centering
\caption{Comparison of iterative schemes for \emph{composite viscoplastic law} \eqref{eq:viscosity-composite}.}\label{tab:composite}
\begin{tabular}{l}
\toprule
Picard solve \textbf{(P)}:\\[1ex]
{$\begin{aligned}
	- \nabla\cdot
	\left[\left(\dfrac{2\viscref\viscstressyield}{2\viscref\strainrateinvII^{\vel_{k-1}}+\viscstressyield}\right)
	\strainratetensor(\bs u_k) \right] + \nabla p_k &= \rhsvel\\
	-\nabla\cdot\vel_k &= 0\\
	\end{aligned}$}\\
\midrule
Standard Newton solve \textbf{(N)}:\\[1ex]
{$\begin{aligned}
	- \nabla\cdot
	\left[\dfrac{2\viscref\viscstressyield}{2\viscref\strainrateinvII^{\vel_{k-1}}+\viscstressyield}
	\bigg(\mathbb I - \dfrac{2\viscref\strainrateinvII^{\vel_{k-1}}}{2\viscref\strainrateinvII^{\vel_{k-1}}+\viscstressyield}
	\dfrac{\strainratetensor(\vel_{k-1})\otimes \strainratetensor(\vel_{k-1})}
	{2(\strainrateinvII^{\vel_{k-1}})^2} \bigg)
	\strainratetensor(\tilde{\vel})\right]
	+ \nabla \tilde p &= -\bs r^{\vel}_{k-1}
	\\
	-\nabla\cdot\tilde{\bs u} &= -r^{p}_{k-1}
	\end{aligned}$}\\
\midrule
Stress--velocity Newton solve \textbf{(SVN)}:\\[1ex]
{$\begin{aligned}
	- \nabla\cdot
	\left[\dfrac{2\viscref\viscstressyield}{2\viscref\strainrateinvII^{\vel_{k-1}}+\viscstressyield}
	\bigg(\mathbb I -
    \dfrac{\left(\strainratetensor(\vel_{k-1})\otimes \viscstresstensor_{k-1}\right)_{\textit{sym}}}
          {2\strainrateinvII^{\vel_{k-1}}\max(\viscstressyield,\viscstressinvII)}
  \bigg)
	\strainratetensor(\tilde{\vel})\right]
	+ \nabla \tilde p &= -\bs r^{\vel}_{k-1}
	\\
	-\nabla\cdot\tilde{\bs u} &= -r^{p}_{k-1}
	\end{aligned}$}\\
Stress update \textbf{(SU)}:\\
{$
	\tilde{\viscstresstensor} =-\viscstresstensor_{k-1} + \dfrac{2\viscref\viscstressyield}{2\viscref\strainrateinvII^{\vel_{k-1}}+\viscstressyield}\left(\strainratetensor(\tilde{\vel}) + \strainratetensor(\vel_{k-1})\right)-
	\dfrac{2\viscref\strainrateinvII^{\vel_{k-1}}}
        {2\viscref\strainrateinvII^{\vel_{k-1}}+\viscstressyield}
	\dfrac{\left(\strainratetensor(\vel_{k-1})\otimes \viscstresstensor_{k-1}\right)_{\textit{sym}}}
        {2(\strainrateinvII^{\vel_{k-1}})^2\max(1,\viscstressinvII/\viscstressyield)}
	\strainratetensor(\tilde{\vel})
$}\\
\bottomrule
\end{tabular}
\end{table}

\section{Test Problems and Implementation}\label{sec:examples}
In this section, we present the test problems we use to study our
algorithms and summarize the discretizations and implementations.
Note that in the test problems described below we use dimensional
units, but our implementation uses rescaled quantities. In particular,
the parameters in Examples 1 and 3 are scaled by $H_0 =$ 30,000 m,
$U_0 = 2.5\times10^{-3}\text{ (m/year)}\times
1/3600/365.25/24 \text{ (year/s) }$ and $\eta_0=10^{22} \text{
  Pa}\cdot\text{s}$. The parameters in Example 2 are scaled by
$H_0 = 1000$ m, $U_0 = 1.0\times10^{-3}\text{ (m/year)}\times
1/3600/365.25/24 \text{ (year/s) }$ and $\eta_0=10^{22}$ Pa\:s.
We use $\bs f=\boldsymbol 0$ in all problems as a gravity force together with
constant density only affects the
pressure but not the velocity for incompressible flow.

\begin{figure}
	\centering
	\includegraphics[width=\columnwidth]{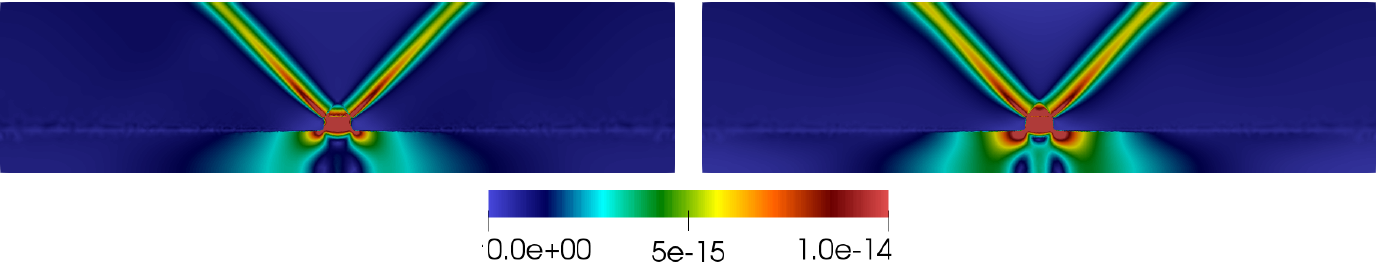}
  \caption{Second invariant of the strain rate pertaining to Example 1 using
           the composite formulation of a von Mises rheology with constant
           yield stress $\tau_y$ (left) and with depth-dependent yield stress
           (right). The depth-dependent yield stress  is given by
           $\tau_y(x,y) = \tau_y\cos(\theta)+\left(\rho g
           y\right)\sin(\theta)$ with
           $\rho = 2700~\text{kg/m}^{3}$, $g = 9.81~\text{m/s}^{2}$, and
           $\theta=30^{\circ}$.
   }
   \label{fig:spiegelman}
\end{figure}

{\bf Example 1: Two-layer compressional notch problem}\quad
The first problem is taken from \citeA{SpiegelmanMayWilson16} and has
also been used in other publications \cite{Kaus10, aspect, FratersBangerthThieulotEtAl19}.
The domain is a rectangle with dimensions $120 \text{ km} \times 30 \text{ km}$.
The fluid consists of a viscoplastic upper layer of thickness $22.5 \text{ km}$ with reference viscosity $\mu_1$ and yield stress $\tau_y$ and an isoviscous lower layer with
constant viscosity $\mu_2$. A notch with the same viscosity as in the
lower layer is introduced in the upper layer. Heterogeneous Dirichlet boundary
conditions in the normal direction are enforced on the left and right
boundaries, $\bs u\cdot \bs{n} = -u_0$, and free-slip conditions are assumed
on the bottom boundary, and thus $\bs u\cdot\bs{n}=0$. All remaining
boundary conditions are homogeneous Neumann conditions. Figure~\ref{fig:spiegelman}
shows the second invariant of the strain rate of solutions with the
composite von Mises rheology. These should be compared with
Fig.~6 in \citeA{SpiegelmanMayWilson16}. We use this problem to study the difference in the solutions computed with ideal and composite rheology,
as well as the convergence behavior of the different algorithms we propose in
this paper.
Unless otherwise specified, we use the parameters  $u_0 = 2.5~\text{mm/yr}$,
    $\mu_1 = 10^{24}~\text{Pa}\:\text{s}$, $\mu_2 =
    10^{21}~\text{Pa}\:\text{s}$, $\viscmin= 10^{19}~\text{Pa}\:\text{s}$ (for ideal rheology) and
    $\tau_y=10^8~\text{Pa}$. The mesh we used for this problem
    (except the runs in Figure \ref{fig:ex1-hp}, where a coarser mesh is used)
    is the same as in
    \citeA{SpiegelmanMayWilson16}. It resolves the notch and consists
    of 36,748 triangles.

\begin{figure}
	\centering
	\includegraphics[width=\columnwidth]{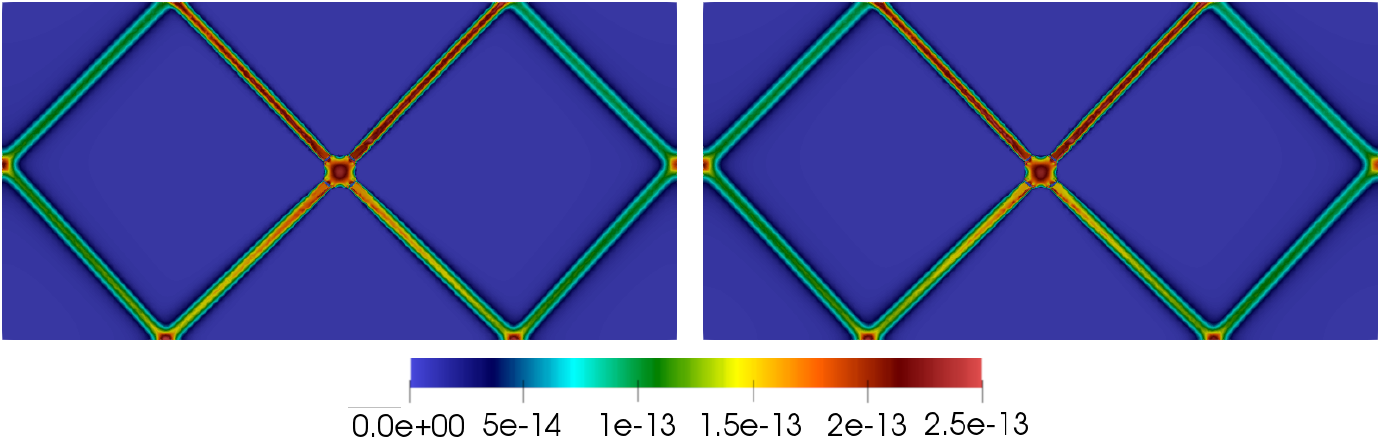}
  \caption{Second invariant of the strain rate pertaining to Example 2 using
           the composite formulation of a von Mises rheology with constant
           yield stress $\tau_y$ (left) and with depth-dependent
           yield stress (right). The depth-dependent
           yield stress is as specified in Figure \ref{fig:spiegelman}.
	}
	\label{fig:duretz}
\end{figure}

\begin{figure}
  \centering
	\includegraphics[width=\columnwidth]{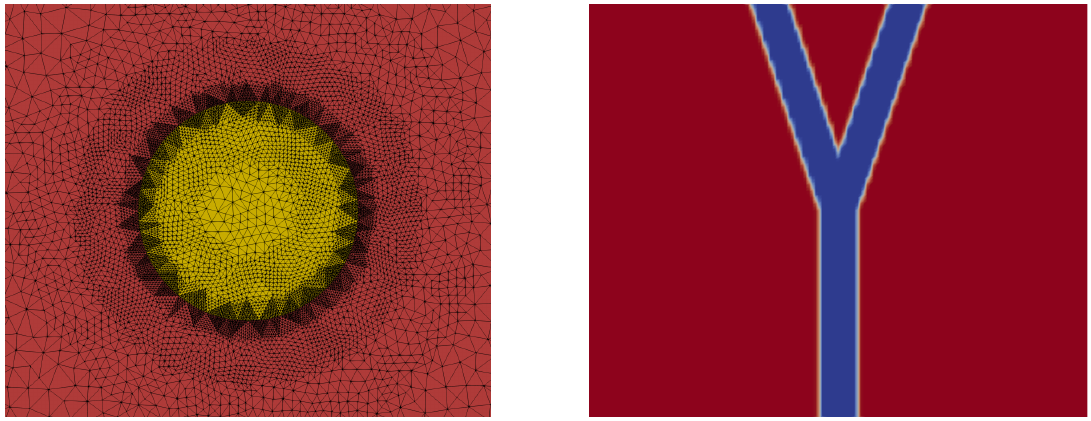}
  \caption{Viscosities of Examples 2 and 3.
    Left: Zoom-in of the locally refined mesh around the
    circular inclusion with radius 0.1 km in Example 2. The viscosity
    outside the the inclusion (in red) is several orders of magnitude
    larger than inside (in yellow).
    Right: Top view of the viscosity of the Y-shaped notch (in blue) of the
    3-dimensional Example 3. Shown is a cutout of size $1\times 1$. The angle
    between the two branches of the notch is $2 \tan^{-1}(1/3) \approx
    37^{\circ}$.
  }
  \label{fig:meshes}
\end{figure}

{\bf Example 2: Circular inclusion problem}\quad
The second test problem is similar to one found in \citeA{DuretzSoucheBorstEtAl18},
where it is used as a test for viscoelastoplastic flow. The domain is a
$4 \text{ km} \times 2 \text{ km}$ rectangle with a circular inclusion (radius
$0.1 \text{ km}$) located in the center. The domain excluding the circular
inclusion has a reference viscosity $\mu_1$ and yield stress $\tau_y$, and
the circular inclusion has a constant viscosity $\mu_2$.
Neumann boundary conditions are applied at the top boundary; Dirichlet
boundary conditions for the normal velocity component are enforced otherwise.
The left and the right boundary conditions are $\bs u\cdot \bs{n} = -u_0$,
resembling lateral compression of the domain, and the bottom is
$\bs u\cdot \bs{n} = u_0/2$.
In Figure \ref{fig:duretz}, we show the second invariant fields
for a problem where the viscosity contrast between the
background and the inclusion is $10^7$. To resolve this extreme
difference, we use a mesh that is significantly refined in and around
the inclusion, as shown on the left in Figure~\ref{fig:meshes}.
We use this setup to study the convergence behavior of our algorithms. For this example, unless otherwise specified, we use parameters $u_0 = 1.0~\text{mm/yr}$, $\mu_1 =
10^{24}~\text{Pa}\:\text{s}$, $\mu_2 =
10^{17}~\text{Pa}\:\text{s}$, $\viscmin = 10^{17}~\text{Pa}\:\text{s}$ (for ideal rheology) and
$\tau_y=3\times10^7~\text{Pa}$.

{\bf Example 3: 3D compressional notch problem}\quad
The third test problem
is a three-dimensional variant of Example 1. The domain is a
$120 \text{ km}\times 60 \text{ km}\times 30 \text{ km}$ rectangular
box, which has a viscoplastic upper layer with reference viscosity
$\mu_1$ and yield stress $\tau_y$ and an isoviscous lower layer with
constant viscosity $\mu_2$. Generalizing the setup from
\citeA{SpiegelmanMayWilson16} used in Example 1, we introduce in the upper layer a Y-shaped notch (see
the top view in Figure~\ref{fig:meshes}, right) with constant
viscosity $\mu_2$. At the left and
right sides, we enforce inflow boundary conditions, $\bs u\cdot \bs{n} = -u_0$,
and at the fore, aft, and bottom boundaries we use $\bs u\cdot \bs{n} = 0$
Dirichlet boundary conditions.  At the top and for tangential velocities we apply
homogeneous Neumann boundary conditions.
With this test problem, we study the
behavior of our algorithms for three-dimensional problems,
their scalability, and their behavior when used in combination with adaptive
mesh refinement of nonconforming hexahedral meshes. The parameters are the
same as in Example 1, namely, $u_0 = 2.5~\text{mm/yr}$, $\mu_1 =
10^{24}~\text{Pa}\:\text{s}$, $\mu_2 = 10^{21}~\text{Pa}\:\text{s}$,
$\viscmin = 10^{19}~\text{Pa}\:\text{s}$,
and $\tau_y=10^8~\text{Pa}$.  Note that in this example,
the edges/faces of mesh elements are not aligned parallel to
the boundary of the notch as they were in the examples before.
Instead, we use a smooth transition from the reference viscosity
$\mu_1$ of the upper layer to the viscosity $\mu_2$ of the notch and lower
layer.  This is carried out by first calculating the (shortest) distance $d$
between a mesh point and the interface and assigning the reference viscosity of
the upper layer to be $\mu_2 + (\mu_1-\mu_2)\exp^{-1000 d^2}$.
This rapidly varying but smooth reference viscosity allows us to represent the physics of
the Stokes problem with a sufficiently refined mesh.  Moreover, along with
solving the nonlinear equations, we adaptively refine the mesh in between
Newton iterations and thus
resolve viscosity gradients that are introduced through the nonlinear rheology.

\begin{figure}
	\centering
	\includegraphics[width=\columnwidth]{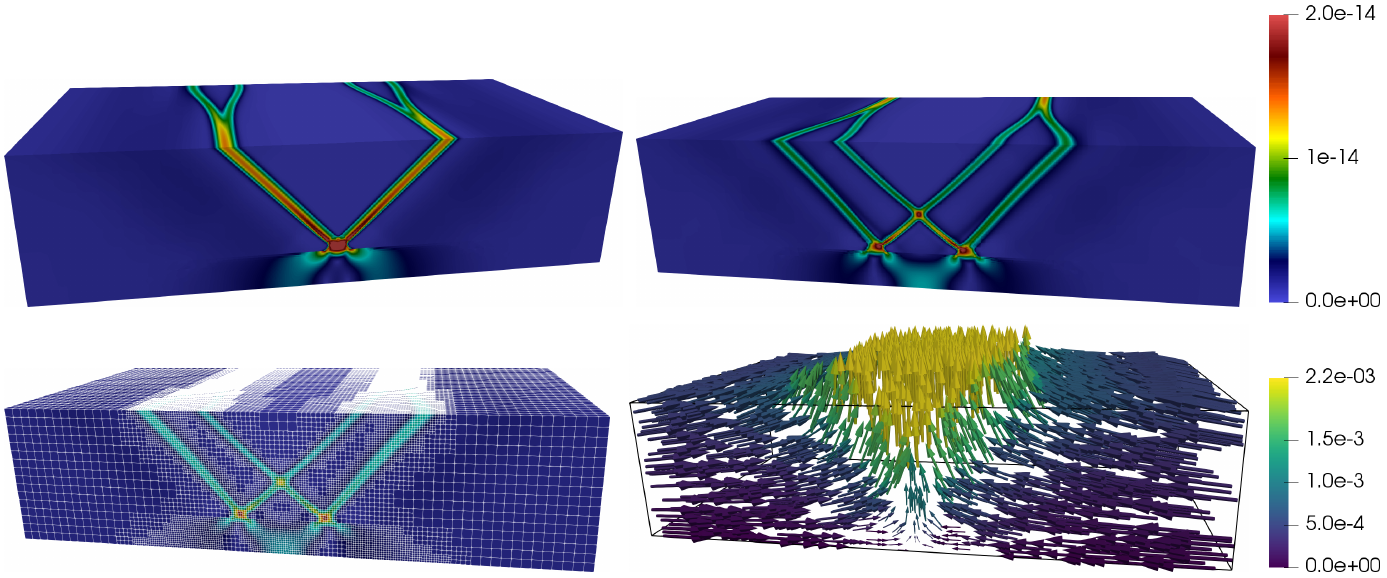}
	\caption{Flow pertaining to Example 3 with ideal viscoplasticity
          rheology. The top row shows the second invariant of the strain rate
          in front view (top left) and back view (top right).
          The bottom row shows the corresponding adaptively refined mesh
          (bottom left) and the velocity field (bottom right) with arrows
          colored by the velocity in $z$-direction.
  }
	\label{fig:Spiegelman3D}
\end{figure}

{\bf Discretization and implementation details}\quad
We now briefly summarize the discretization and implementation we
use for our computational results. While the development of efficient
and scalable solvers and preconditioners for (linearized) Stokes
problems is an interesting and active field of research, here we focus
on the behavior of the nonlinear iterations and give only a brief
summary of the computational methods used.  We discretize the Stokes
equations using stable finite element pairings for the velocity and
pressure variables.  Note that one could use different (stable or
stabilized) finite element pairings from those used here, or even different
discretization methods such as staggered finite differences
\cite{DabrowskiKrotkiewskiSchmid18}.

For the two-dimensional problems, unless stated otherwise, we use
standard $P_2\times P_1$ Taylor-Hood elements on triangular meshes
\cite{ElmanSilvesterWathen14a} and discontinuous linear elements for
the stress variable. We also report results with higher-order Taylor
Hood elements. Additionally, we experimented with stable element
pairings with discontinuous pressure spaces. Since we observed
very similar convergence behavior, we decided to not show these
results. Our two-dimensional implementation
\footnote{Source code available at \url{https://bitbucket.org/johannrudi/perturbed_newton}.}
is based on the open source finite element library FEniCS
\cite{LoggMardalGarth12} and thus uses weak forms and numerical
quadrature to compute finite element matrices. We found it beneficial
to the convergence of iterative methods (in particular,
Newton-type methods) to increase the numerical quadrature order, for example to
10. Doing so results in small quadrature errors and thus accurately
approximated weak forms. These are important when linearization is done
on the weak form level (as in FEniCS), to avoid having linearization and
quadrature interfere and small numerical errors spoil the fast local convergence
of Newton-type methods. Even with lower quadrature, however,
we find good convergence but usually observe suboptimal convergence
close to the solution. The linear problems arising in each Newton step
are solved with a direct solver. To compute an appropriate step length
for the Newton updates, we check whether the optimization objective in
\eqref{eq:minJ} decreases, and we reduce the step length if this is not the case.
We compared using the optimization objective with using the nonlinear
residual in this step length criterion and found that the two
methods behave similarly.

For the three-dimensional problem, we use $\mathbb Q_2\times
P_{1}^{\text{disc}}$ elements on nonconforming hexahedral meshes
\cite{ElmanSilvesterWathen14a, RudiStadlerGhattas17}. The dual stress
variable is treated pointwise at quadrature nodes and not
approximated by using finite elements. This implementation builds on our
previously developed scalable parallel software framework
\cite{RudiStadlerGhattas17, RudiMalossiIsaacEtAl15}. To solve the
(large) linear systems that arise, we use an iterative Krylov method with
a block preconditioner for the Stokes system, which entails
a Schur complement preconditioner.  The viscous block in the Stokes system and
components in our Schur complement preconditioner are amenable to
multigrid methods.  Hence, we have developed a hybrid (geometric
and algebraic) multigrid method for adaptive meshes \cite{RudiStadlerGhattas17,
RudiMalossiIsaacEtAl15}. For these three-dimensional problems,
the step length for the Newton update is found by checking
descent in the $L^2$-norm of the nonlinear residual, which amounts to
approximately inverting a mass matrix during the inner product:
  $\norm{\boldvar{r}}_{L^2} = \norm{\vec{r}}_{\mat{M}^{-1}} =
   \sqrt{\vec{r}\transpose\mat{M}^{-1}\vec{r}}$,
for a residual vector $\vec{r}$ and the mass matrix $\mat{M}$ of the
corresponding finite element space.
In practice, it is beneficial for time to solution to employ an inexact
Newton--Krylov method: far from the nonlinear solution, the stopping
tolerance of the Krylov linear solver is set to a relatively high value,
whereas closer to the nonlinear solution, the tolerance is decreased.  Since
this paper focuses on Newton's method, we set a low Krylov tolerance $<10^{-7}$
throughout the nonlinear solve in our numerical examples.

The next section discusses qualitative properties of the solutions to the test
problems and compares the convergence of the algorithms from
Section~\ref{sec:Newton}.

\section{Numerical Results}\label{sec:results}
In this section, we study the performance of our algorithms for the two-
and three-dimensional test problems described above. We also study the differences in the
solutions resulting from using the ideal and the composite viscoplastic
rheology laws, and we numerically verify the energy minimization properties of
Stokes solutions derived in the preceding sections.

\subsection{Comparison between ideal and composite plasticity}
\begin{figure}
  \centering
	\includegraphics[width=\columnwidth]{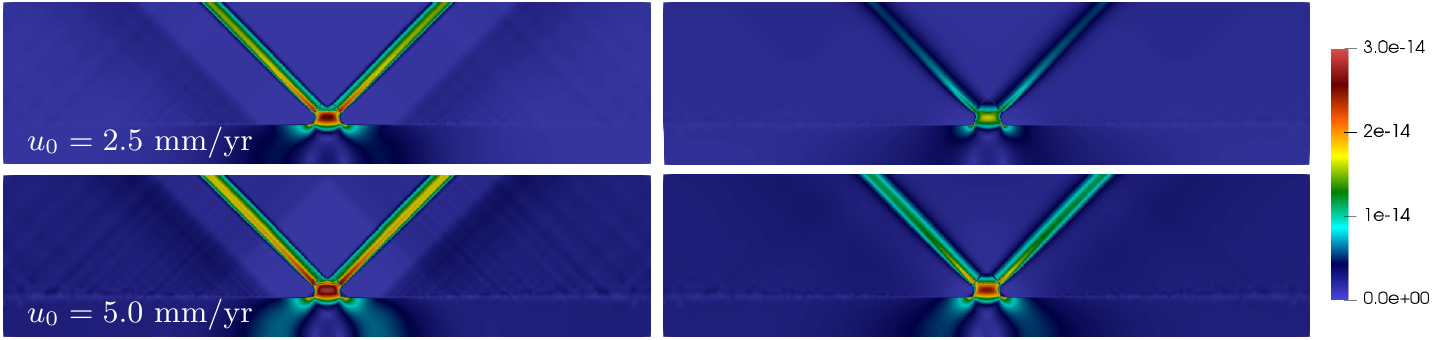}
\caption{Comparison between ideal (left column) and composite rheology (right
  column) for von Mises rheology with constant yield stress and for two
  different values of lateral inflow velocities $u_0$.
  Colors depict the second invariant of the strain
  rate for Example 1.\label{fig:Spiegelman-compare}}
\end{figure}

\begin{figure}
  \centering
	\includegraphics[width=\columnwidth]{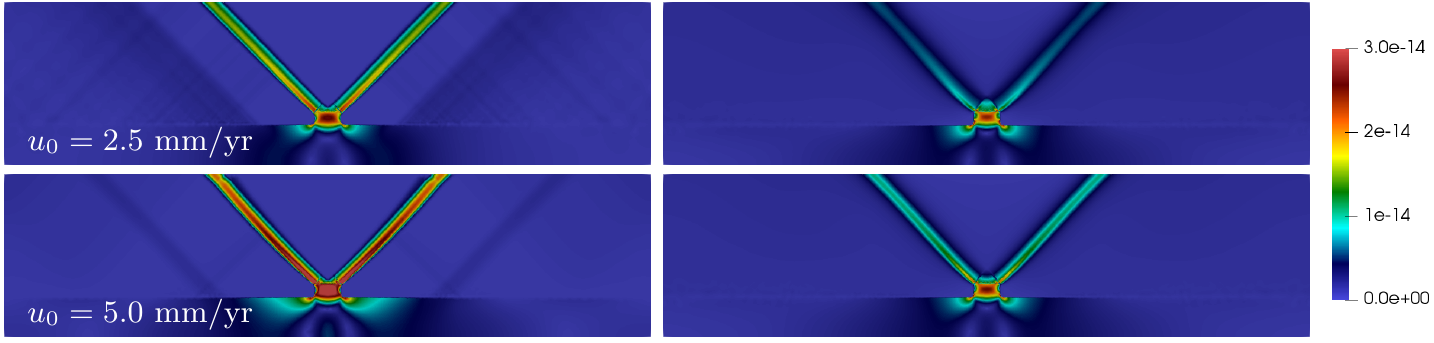}
\caption{Comparison between ideal (left column) and composite rheology
  (right column) for depth-dependent von Mises rheology and for two different
  values of lateral inflow velocities $u_0$. The values for the depth-dependent
  yield stress are as specified in Figure \ref{fig:spiegelman}.
  Colors depict the second
  invariant of the strain rate for Example 1.\label{fig:Spiegelman-compare-depth}}
\end{figure}

First, we compare the solution structure obtained with
\eqref{eq:viscosity-ideal} and \eqref{eq:viscosity-composite}, the
ideal and the composite von Mises rheology laws, respectively.  Figures
\ref{fig:Spiegelman-compare} and \ref{fig:Spiegelman-compare-depth}
show the second invariants of the strain rates for Example 1 when
solved by using these two variants of combining the yield and
reference viscosities.
As illustrated in Figure~\ref{fig:Phi} and as discussed previously in
in Appendix C in \citeA{SpiegelmanMayWilson16},
the composite rheology law can be considered as a
regularization of the ideal plasticity due to a smooth transition
between yielding and reference viscosity. Since the composite
viscosity is always below the ideal viscosity, it leads to smaller
strain rates in the solutions, as demonstrated in Figures
\ref{fig:Spiegelman-compare} and
\ref{fig:Spiegelman-compare-depth}. This difference in magnitude
appears to be stronger for the problem with depth-dependent von Mises
rheology (Figure~\ref{fig:Spiegelman-compare-depth}). Additionally, the shape
of the fault zone differs slightly, in particular close to the notch;
and the second invariant of the strain rate is smoother for
the composite rheology. This is, again, likely a consequence of the
smooth way the yield viscosity is combined with the reference
viscosity in the composite rheology law. This comparison shows that the
definition of the viscoplastic rheology may impact the magnitude and
potentially the structure of the solution; thus the two models cannot be
interchanged.

\begin{figure}
  \centering
	\includegraphics[width=\columnwidth]{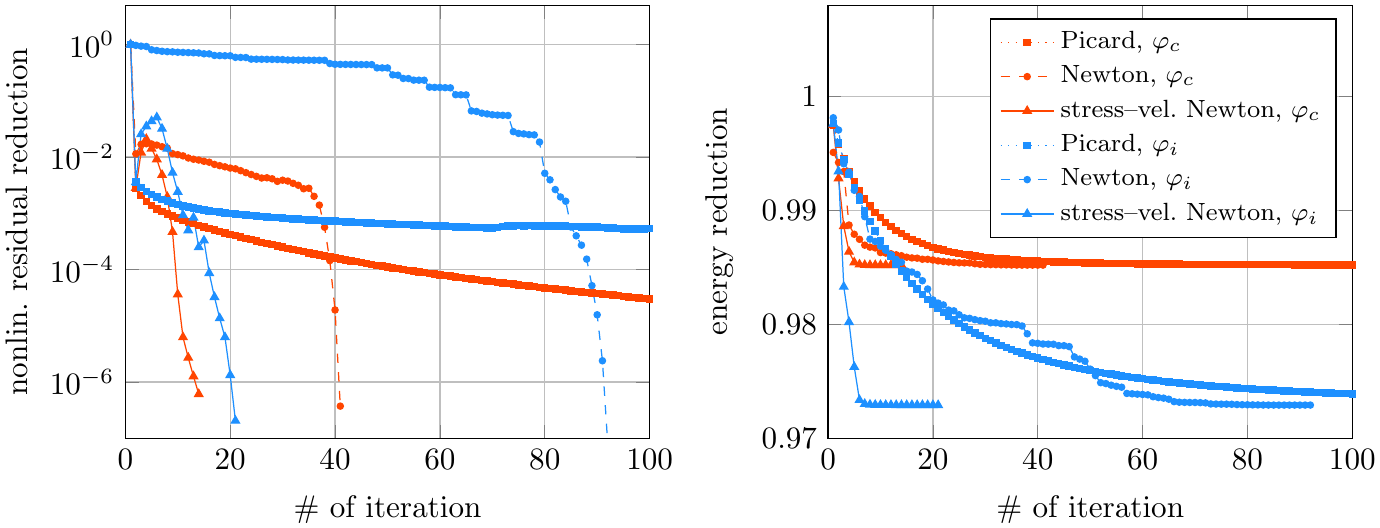}
  \caption{Convergence behavior of Example 1 for different iterative
    methods and ideal (blue curves) and composite (red curves) viscoplasticity.
    Shown on the left is the nonlinear residual reduction in each iteration
    and on the right the corresponding reduction of the energy objective
    \eqref{eq:minJ}. }
  \label{fig:conv-ex1}
\end{figure}

\subsection{Results for the three-dimensional problem}
\label{sec:results-ex3}
Numerical solutions for Example 3, a three-dimensional version of the notch
problem, are shown in Figure \ref{fig:Spiegelman3D}.
The symmetric V-shaped faults from Example 1 split
into thinner W-shaped faults due to the shape of the notch
along the interface of
the upper and lower layers (see Figure~\ref{fig:meshes}, right).
Figure \ref{fig:Spiegelman3D}, bottom right shows the velocity field (green
arrows) at the nonlinear solution, where the central region is uniformly
lifted with two lateral secondary uplift regions due to the secondary faults
caused by the geometry of the notch.

We performed the three-dimensional computations on two uniformly refined
meshes, including a relatively coarse mesh, which sufficiently resolves
physical features, of 262,144 elements and a fine mesh of 2,097,152 elements,
which is generated by one additional level of refinement.
Further, we performed nonlinear solves with adaptive mesh refinement in
between Newton iterations.
The mesh adaptation starts from a
uniform mesh of 262,144 elements and coarsens and refines based on the average
second invariant of the strain rate per element:
  $\frac{1}{|\Omega_e|}\int_{\Omega_e} \strainrateinvII^{\bs u}(\bs x) d\bs x$,
where $|\Omega_e|$ is the volume of an element $e$.
The adaptation is triggered after a Newton step until the fifth nonlinear
iteration, if it is needed.  The final adapted meshes have a maximum refinement
that is one level finer than the fine uniform mesh mentioned before.  They
differ in size for standard and stress--velocity Newton due to differences in
the convergence of the solvers, having 881,945 and 1,107,758 elements,
respectively.
One example of an adapted mesh is depicted by the white edges in Figure
\ref{fig:Spiegelman3D}, bottom left.

\subsection{Nonlinear solver convergence}\label{sec:solver-conv}
\begin{figure}
  \centering
	\includegraphics[width=\columnwidth]{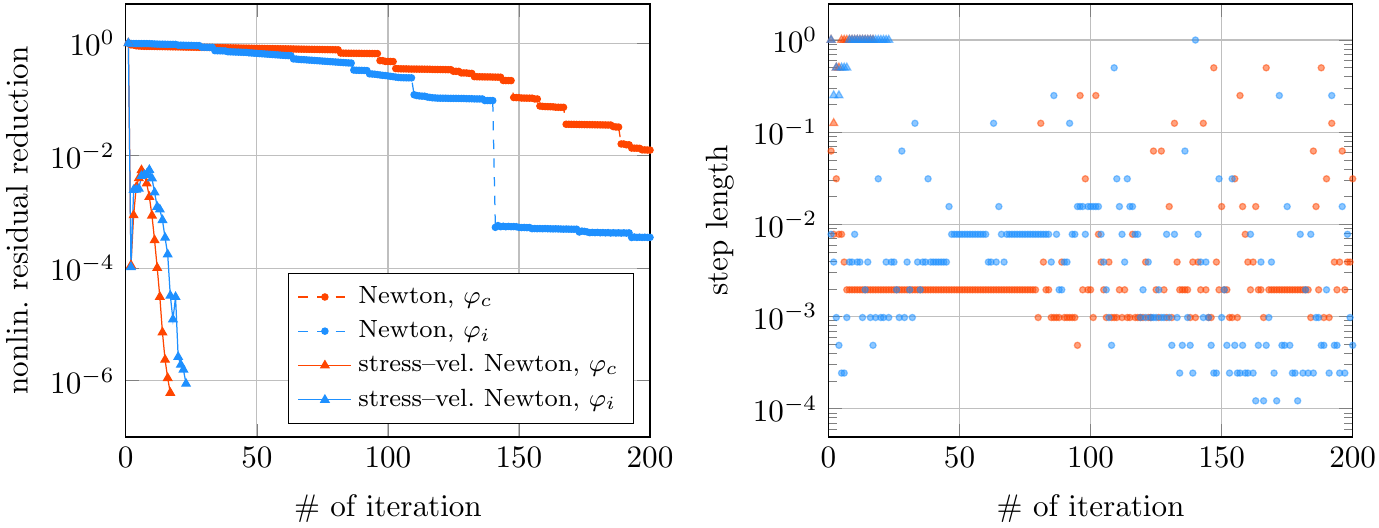}
  \caption{Convergence behavior for Example 2. Shown on the left is the residual
           reduction for standard and stress--velocity Newton for ideal (blue
           curves) and composite (red curves) rheologies. For ideal
           viscoplasticity, the lower viscosity bound $\viscmin$ is important
           for convergence of the standard Newton method. If it is reduced by a
           factor of 10, the line search fails after 65 iterations (the
           stress--velocity Newton method still converges).
           The right plot shows the update step sizes these methods take.
  }
  \label{fig:conv-ex2}
\end{figure}

\begin{figure}
  \centering
	\includegraphics[width=0.9\columnwidth]{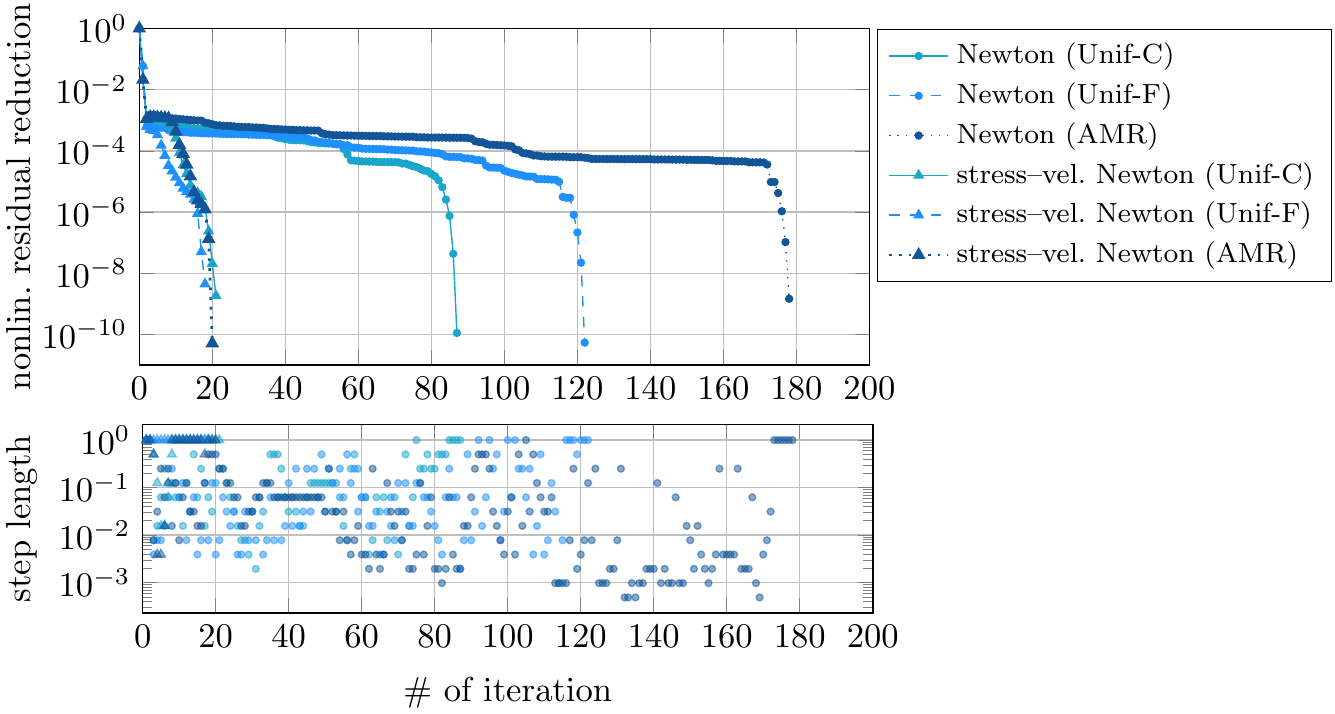}
  \caption{Convergence results for the three-dimensional Example 3 using the
           ideal viscoplastic rheology $\varphi_i$. The coarse uniform mesh
           (Unif-C) consists of 262K %
           elements, the fine uniform mesh (Unif-F) with one additional level
           of refinement has 2.1M %
           elements, and the final adapted meshes (AMR) with another additional
           level of max refinement have $\sim$1M elements.
           We observe a mesh-independent convergence of stress--vel.\ Newton.
           It is also faster than standard Newton by a factor of 4--9, where
           convergence deteriorates with finer mesh resolutions.
  }
	\label{fig:conv-ex3}
\end{figure}

Next
we study the convergence behavior of the different nonlinear iterative methods.
Figures \ref{fig:conv-ex1}, \ref{fig:conv-ex2}, and \ref{fig:conv-ex3}
show the convergence history for Examples 1, 2, and 3.
First, we observe that
all algorithms converge faster for the composite compared with the ideal rheology
(red vs.\ blue curves in Figure \ref{fig:conv-ex1}). This result is likely due to the
higher regularity of the composite
rheology \eqref{eq:viscosity-composite} compared with the ideal rheology
\eqref{eq:viscosity-ideal}, which involves the pointwise
$\min$-function.  As is well documented in the literature, the Picard
fixed-point method initially reduces the nonlinear residual rapidly
but makes very slow progress in later iterations (see Figure
\ref{fig:conv-ex1}), in particular for the highly nonlinear ideal rheology. The
standard Newton method makes very slow progress initially until
it reaches the ``basin'' of fast convergence, where it performs with the expected
fast convergence. Furthermore the figure demonstrates that the stress--velocity
Newton method converges rapidly for both types of viscoplastic laws (a factor
of 2.5--5 fewer iterations than standard Newton), which is one of the key
advantages of these novel nonlinear solvers.
Note that the numerical results support the theoretical observation at the end
of Section \ref{sec:pdNewton} that, if the stress variable is initialized
with zero, the first iteration of stress--velocity Newton coincides with the
first iteration of the Picard method.

The main reason for the vastly different behavior between the standard Newton
and the stress--velocity Newton method can be seen in the step length plot of
Figures \ref{fig:conv-ex2}, right.  The line search reduces the
length of standard Newton steps drastically in most iterations except the the
very end, while the stress--velocity Newton method mostly makes full steps of
length 1. A heuristic explanation for this behavior is that by introducing the
additional stress variable into the nonlinear Stokes problem and using the
transformation \eqref{eq:NCPs}, the system becomes ``less'' nonlinear.

The results for Example 1 also show that while the
nonlinear residual is not reduced monotonically, the energy objective
\eqref{eq:minJ} decreases monotonically (Figure \ref{fig:conv-ex1}, right).
Since the objective is used
to check descent when using the Newton update, this shows that using
the objective from the underlying optimization problem for descent
might be superior to using the nonlinear residual because it allows larger
steps.

Our observations of the two-dimensional examples carry over to the
three-di\-men\-sio\-nal example in Figure \ref{fig:conv-ex3}.  The top graph of
this figure shows the nonlinear convergence using the computationally more
challenging ideal viscoplastic rheology.
We have performed numerical experiments on three mesh configurations with
increasing resolution comprised of a relatively coarse uniform mesh (Unif-C), a
fine uniform mesh (Unif-F), and adaptively refined meshes (see Section
\ref{sec:results-ex3} for more details).  The finest mesh is achieved with the
adaptively refined setup. It has an additional level of refinement to
(Unif-F) that is assumed only locally resulting in fewer mesh elements overall
than with uniform refinement.  The (Unif-C) mesh consists of 262K elements, the
(Unif-F) mesh with one additional level of refinement has 2.1M elements, and
the (AMR) meshes
have $\sim$1M elements.
The stress--velocity Newton method takes about 20 iterations to reduce the
nonlinear residual by 8--10 orders of magnitude for all three mesh
configurations.  The corresponding step lengths of stress--velocity Newton
reduce below one only during the initial eight iterations (i.e., far from the
solution) due to backtracking line search, but they remain at one until
the algorithm has converged.
The convergence of standard Newton is slower by a factor of four for the
(Unif-C) mesh and increases to a factor of nine for the (AMR) setup, which has
the highest level of local mesh refinement.  Simultaneously with stagnating
convergence, we observe a dramatic reduction of step lengths until a somewhat
random inflection point when the desired super-linear Newton convergence
reduces the residual quickly.

\begin{table}
  \caption{Convergence with varying mesh size for the composite rheology of
           Example 2.  We list the number of iterations to reduce the nonlinear
           residual by $10^{-6}$.
  }
  \label{tab:conv-ex2}
  \centering
	\begin{tabular}{lcccc}
		\toprule
		\# mesh elements                   &2,140 & 7,416 & 24,617 & 91,350\\
    \midrule
    \# standard Newton iterations      &53 &          104 &          287 &          400\\
    \# stress--vel.\ Newton iterations &18 &\phantom{1}19 &\phantom{2}19 &\phantom{4}17\\
    \bottomrule
	\end{tabular}
\end{table}

\subsection{Solver behavior under mesh refinement and increasing order}
\begin{figure}
	\centering
	\includegraphics[width=\columnwidth]{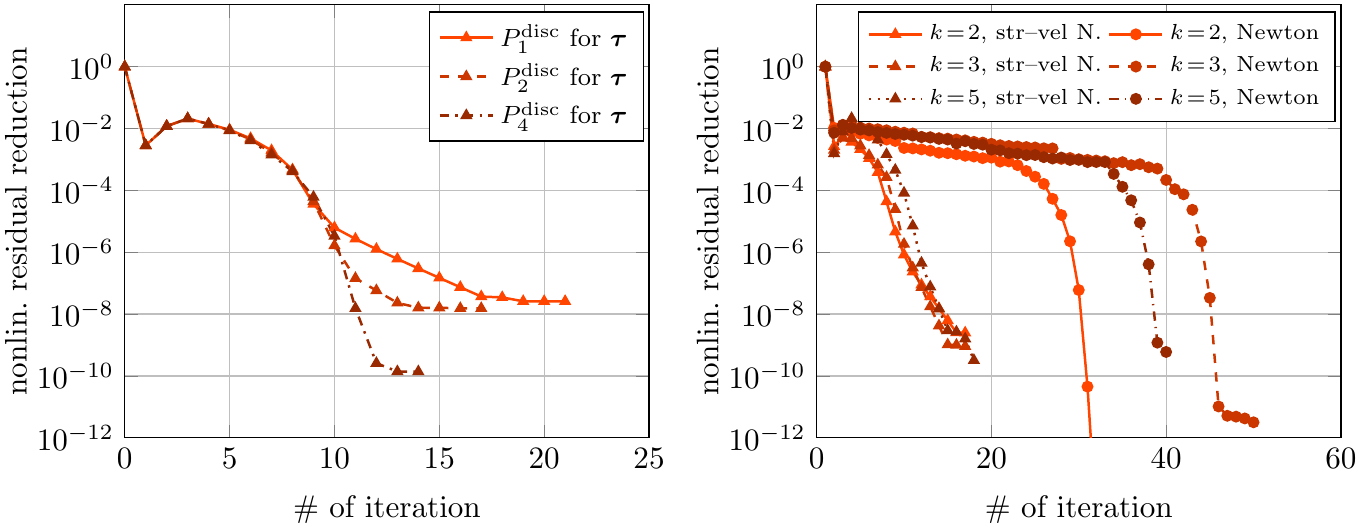}
        \caption{Tests with higher-order discretizations for Example 1
          with composite rheology and parameters as detailed in Figure
          \ref{fig:conv-ex1}. Shown on the left is the residual
          reduction using standard Taylor-Hood elements, namely, order
          $k=2$ for the velocity and $k-1$ for the pressure, but with
          different-order discontinuous elements to discretize the
          stress tensor $\viscstresstensor$.
          Shown
          on the right is the behavior of the standard and the stress--velocity Newton
          methods for higher-order Taylor-Hood discretizations, namely,
          order $k$ for velocity and order $k-1$ for the pressure. The
          stress tensor is always discretized with polynomials of degree
          $k-1$. These runs use a coarser mesh consisting of 1,892 elements.
        }\label{fig:ex1-hp}
\end{figure}

In particular for solving large-scale problems, it is crucial that
iterations numbers are insensitive to the refinement of the mesh and
also the polynomial degree. Thus, in Table \ref{tab:conv-ex2} we study the
behavior of the standard and
stress--velocity Newton method for meshes with an increasing number of
elements. As can be seen, the iteration number for the standard Newton
method grows significantly, while the iteration number for the
velocity--stress Newton method remains constant. A similar test, now for
increasing the order of the discretizations while keeping the mesh
fixed, is shown on the right of Figure \ref{fig:ex1-hp}. We find that
the convergence behavior for the stress--velocity Newton method does
not change significantly when the polynomial order for the velocity
discretization is increased from order $k=2$ to $k=5$, while the
number of iterations increases with the polynomial order for the
standard Newton method.

Supporting our results for two-dimensional problems using the FEniCS library,
our code for three-dimensional problems shows similar convergence properties
with respect to mesh refinement in Figure \ref{fig:conv-ex3}.
It shows three mesh configurations with increasing level of refinement, where
(Unif-C) and (Unif-F) are coarse and fine uniform meshes, respectively, and
(AMR) has one additional level of local adaptive refinement compared to
(Unif-F) (see Sections \ref{sec:results-ex3} and \ref{sec:solver-conv} for more
details).
Stress--velocity Newton requires only about 20 iterations until convergence and
remains independent of mesh size.  In contrast, the iteration counts of
standard Newton depend heavily on the mesh resolution and increase by factor of
$\sim$1.5 each time the mesh gains an additional level of refinement:
taking 87, 122, and 178 iterations for (Unif-C), (Unif-F), and (AMR) meshes,
respectively.
These numerical results demonstrate the significant advantages of the
stress--velocity Newton method for large-scale three-dimensional viscoplastic
problems.

\section{Conclusions}
While viscoplastic rheologies are important for numerical models of the
lithosphere, the resulting nonlinearity in the equations
has been challenging for computational solutions.
For viscoplastic rheologies where the yield stress
does not depend on the dynamic pressure, we propose a novel
Newton-type method in which we formally introduce an independent
stress variable, compute the Newton linearization, and then use a Schur
complement argument to eliminate the update for the stress variable
from the Newton system. Hence, our method requires solving a
sequence of linear Stokes problems that have the usual
velocity--pressure unknowns. This method mostly allows one to take full
unit length Newton update steps and thus converges significantly
faster than the standard Newton method and the Picard method; and it is at least as
fast as combined Picard--Newton methods, which require manual tuning and do not generalize
well to different problem setups. We also observe
that the number of iterations the method requires is independent of how fine
the discretization is and which polynomial order is used to approximate
velocity and pressure.

We also compare two von Mises rheology formulations: an ideal variant
that uses a pointwise $\min$-function and a composite variant that
uses the harmonic mean and, therefore, has effective viscosities
that always lie below the ideal variant. While the two variants are sometimes
used interchangeably, we find that they lead to solutions of rather
different magnitude. Because of the improved smoothness, all solution
algorithms in our comparison converge faster for the composite than for the
ideal rheology. Furthermore, we study the uniqueness of solutions pertaining to the two
rheologies using variational energy minimization arguments. We find
that problems with the composite rheology always have a unique
solution, whereas problems with ideal von Mises viscoplasticity require
a regularization to ensure that they have a unique solution. In
Example 2, we find that this regularization parameter can also play a
practical role for the convergence of the standard Newton algorithm;
see Figure~\ref{fig:conv-ex2}.

Systematic study is needed to determine whether
the novel methods
can be extended to
obtaining better converged solutions for rheologies with
generalizations of the von Mises viscoplasticity, such as
viscoelastoplasticity and rheologies where the yield stress depends on
the dynamic pressure, also referred to as Drucker--Prager.

\section*{Acknowledgments}
We thank Anton Popov and Dave May for insightful reviews and the resulting
improvements of this paper.

This material is based upon work supported by the U.S.\ Department of Energy,
Office of Science, Advanced Scientific Computing Research under Contract
DE-AC02-06CH11357 and the Exascale Computing Project (Contract No.\
17-SC-20-SC).  This research was supported by the Exascale Computing Project
(17-SC-20-SC), a collaborative effort of the U.S.\ Department of Energy Office
of Science and the National Nuclear Security Administration.

The work was partially supported by the US National Science Foundation through
grant EAR-1646337 and DMS-1723211. Computing time on TACC's Stampede2
supercomputer was provided through the Extreme Science and Engineering
Discovery Environment (XSEDE), which is supported by National Science
Foundation grant number ACI-1548562.

The source code of our two-dimensional implementation, which is based on the
finite element library FEniCS, is available in the git repository:
\url{https://bitbucket.org/johannrudi/perturbed_newton}

\bibliographystyle{plain}
\bibliography{ccgo,extra}

\newpage
\noindent
{\bf Government License}\quad
The submitted manuscript has been created by UChicago Argonne, LLC, Operator of
Argonne National Laboratory (``Argonne''). Argonne, a U.S.\ Department of Energy
Office of Science laboratory, is operated under Contract No.\ DE-AC02-06CH11357.
The U.S.\ Government retains for itself, and others acting on its behalf, a
paid-up nonexclusive, irrevocable worldwide license in said article to
reproduce, prepare derivative works, distribute copies to the public, and
perform publicly and display publicly, by or on behalf of the Government.  The
Department of Energy will provide public access to these results of federally
sponsored research in accordance with the DOE Public Access Plan.\\
\url{https://energy.gov/downloads/doe-public-access-plan}

\end{document}